\magnification=\magstep1
\vsize=23.5true cm
\font\b=cmssdc10 scaled\magstep1
\nopagenumbers
\topskip=1true cm
\headline={\tenrm\hfil\folio\hfil}
\raggedbottom
\abovedisplayskip=3mm
\belowdisplayskip=3mm
\abovedisplayshortskip=0mm
\belowdisplayshortskip=2mm
\normalbaselineskip=12pt
\normalbaselines
\parindent=0cm
\def\Wedge#1{\buildrel\wedge\over#1}
\def\Vee#1{\buildrel\vee\over#1}
\def\Sim#1{\buildrel\sim\over#1}
\def\sqr#1#2{{\vcenter{\vbox{\hrule height.#2pt
\hbox{\vrule width.#2pt height#1pt \kern#1pt
\vrule width.#2pt}
\hrule height.#2pt}}}}
\def\square{\mathchoice\sqr34\sqr34\sqr{2.1}3\sqr{1.5}3}
%NOTHING PERMITTED BETWEEN HERE AND TEXT!

%TEXT STARTS HERE
%REVISED 2-3-1997
\centerline{\bf An Uncountable Family of Regular Borel Measures}
\centerline{\bf On Certain Path Spaces of Lipschitz Functions.}
\bigskip
\centerline{01/19/07}
\bigskip
\bigskip
\centerline{\it R. L. Baker}
\centerline{\it University of Iowa}
\centerline{\it Iowa City, Iowa 52242}
\vfill\eject
\centerline{\bf An Uncountable Family of Regular Borel measures} 
\bigskip
\centerline{\it R. L. Baker}
\centerline{\it  University of Iowa}
\centerline{\it  Iowa City, Iowa 52242}
\vfill\eject
\noindent
ABSTRACT. Let $c>0$ be a fixed constant. Let $0\le r<s$ be an arbitrary pair of real numbers. Let $a,b$ be any pair of real numbers such that
$|\,b-a\,|\le c(s-r)$. Define $C^s_r$ to be the set of continuous real-valued functions on $[r,s]$, and 
define $C_r$ to be the set of continuous real-valued functions
on $[\,r,+\infty)$. Finally, consider the following sets of Lipschitz functions:
$$\eqalignno{\Lambda^s_r &=
                \{\,x\in C^s_r\;|\;|x(v)-x(u)|\le c|v-u|,\quad\hbox{\sl for all}\quad u,v\in[\,r,s\,]\,\},&(1)\cr
    \Lambda_r &=\{\,x\in C_r\;|\;|x(v)-x(u)|\le c|v-u|,\quad\hbox{\sl for all}\quad u,v\in[\,r,+\infty)\,\},&(2)\cr
    \Lambda^{s,b}_{r,a}&=\{\,x\in\Lambda^s_r\;|\;x(r)=a,\quad x(s)=b\,\},&(3)\cr
    \Lambda^s_{r,a}&=\{\,x\in\Lambda^s_r\;|\;x(r)=a\,\},&(4)\cr
    \Lambda^{s,b}_r&=\{\,x\in\Lambda^s_r\;|\;x(s)=b\,\},&(5)\cr
    \Lambda_{r,a}&=\{\,x\in\Lambda_r\;|\;x(r)=a\,\}.&(6)\cr}$$
We present a general method of constructing an uncountable family of regular Borel measures on each of the sets (1), (2),
and an uncountable family of regular Borel probability measures on each of the sets (3)-(6).
Using this method, we give a definition of {\bf Lebesgue measure} on the sets (1) and (2), and a definition of
{\bf the uniform probability} measure on each of the sets (3)-(6). 

\noindent
{\bf Mathematical Reviews subject classification:} 26A99, 28C05, 28C15, 28C20, 60G05, 81S40.
\vfill\eject
\parindent=1cm
\centerline{\bf 1. INTRODUCTION}
\vskip 1cm
\noindent
Let $c>0$ be a fixed constant. Let $\hbox{\rm\b R}$ be the set of real numbers. For each interval 
$I\subseteq\hbox{\rm\b R}$, let $\Lambda(I)$ be the set of functions $x:I\to\hbox{\rm\b R}$ such that
$x$ satisfies the Lipschitz condition
$$|x(t)-x(s)|\le c|t-s|,\quad \hbox{\rm for all},\quad s,t\in I.$$
Let $r,s$ be any pair of real numbers such that $0\le r<s$. Define

$$\Lambda^s_r=\Lambda([r,s]),\quad \Lambda_r=\Lambda([r,+\infty)).$$
Finally, let $a,b$ be any pair of real numbers such that $|b-a|\le c(s-r)$. Then define
$$\eqalignno{\Lambda^{s,b}_{r,a}&=\{\,x\in\Lambda^s_r\;|\;x(r)=a,\quad x(s)=b\,\};&(1)\cr
             \Lambda^s_{r,a}&=\{\,x\in\Lambda^s_r\;|\;x(r)=a\,\};\cr
	     \Lambda^{s,b}_r&=\{\,x\in\Lambda^s_r\;|\;x(s)=b\,\};\cr
	     \Lambda_{r,a}&=\{\,x\in\Lambda_r\;|\;x(r)=a\,\}.\cr}$$

The main result of the present paper (Sections 2 and 3) is a general method of constructing an uncountable 
family of regular Borel measures 
$$\displaylines{\hfill \lambda^{s,b}_{r,a},\quad \lambda^s_{r,a},\quad \lambda^{s,b}_r,\quad \lambda_{r,a},
                \hfill\llap{\rm (2)}\cr
                \hfill \lambda^s_r,\quad \lambda_r,
		\hfill\llap{\rm (3)}\cr}$$
respectively, on each of the following sets:
$$\displaylines{\hfill \Lambda^{s,b}_{r,a},\quad  \Lambda^s_{r,a},\quad  \Lambda^{s,b}_r,\quad \Lambda_{r,a};
                 \hfill\llap{\rm (4)}\cr
		  \hfill \Lambda^s_r,\quad \Lambda_r.  \hfill\llap{\rm (5)}\cr}$$
Each Borel measure in (2) is a probability measure, constructed as the continuous image of Lebesgue measure under a certain
family of continuous surjective mappings
$$\varphi^{s,b}_{r,a}:\,\Omega^s_r\to\Lambda^{s,b}_{r,a},\quad
  \varphi^s_{r,a}:\,\Wedge\Omega\mathstrut^s_r\to\Lambda^s_{r,a},\quad
  \varphi^{s,b}_r:\,\Vee\Omega\mathstrut^s_r\to\Lambda^{s,b}_r,\quad
  \Wedge\Omega_r\to\Lambda_{r,a}. \eqno (6)$$
Each of the spaces
$$\Omega^s_r,\quad \Wedge\Omega\mathstrut^s_r,\quad \Vee\Omega\mathstrut^s_r,\quad \Wedge\Omega_r \eqno (7)$$
is endowed with Lebesgue measure, and has the form $[0,1]^A$, where $A$ is some indexing set.
Likewise, each Borel measure in (3) is constructed as the continuous image of Lebesgue measure under a certain family of
continuous surjective mappings
$$\varphi^s_r:\,\Sim\Omega\mathstrut^s_r\to\Lambda^s_r,\quad \varphi_r:\,\Sim\Omega_r\to\Lambda_r, \eqno (8)$$
where each of the spaces
$$\Sim\Omega\mathstrut^s_r,\quad \Sim\Omega_r \eqno (9)$$
has the form $[0,1]^B$ for some indexing set $B$, and is endowed with Lebesgue measure.

In Section 4, certain members 
$$\lambda^{s,b}_{r,a},\quad \lambda^s_{r,a},\quad \lambda^{s,b}_r,\quad \lambda_{r,a}$$
of the uncountable family (2) are singled out and defined to be {\bf the uniform probability measure} on the spaces
$$\Lambda^{s,b}_{r,a},\quad  \Lambda^s_{r,a},\quad  \Lambda^{s,b}_r,\quad \Lambda_{r,a}, \eqno (10)$$
and certain members
$$\lambda^s_r,\quad \lambda_r$$
of the uncountable family (3) are singled out and defined to be {\bf Lebesgue measure} on the spaces
$$\Lambda^s_r,\quad \Lambda_r. \eqno (11)$$

\vskip 1cm
\centerline{\bf 2. CONSTRUCTION OF THE FUNCTIONS $\varphi^{\,\square}_{\,\square}$}
\vskip 1cm
\noindent
In this section we construct the families of continuous surjective mappings mentioned in (6) and 
(1) of the introduction.
\medskip
\proclaim Definition 2.1. Let $0\le r<s$ be given. Let $(r,a),(s,b)$ be two given points in the plane. Define
$$\eqalign{ F_{r,a}&=\{\,(t,x)\,|\,r\le t \quad {\rm and} \quad  a-c(t-r)\le x\le a+c(t-r)\,\},\cr
            B_{s,b}&=\{\,(t,x)\,|\,t\le s \quad {\rm and } \quad b-c(s-t)\le x\le b+c(s-t)\,\},\cr
	    P^{s,b}_{r,a}&=F_{r,a}\cap B_{s,b}.\cr}$$
\medskip

\noindent
{\bf Proposition 2.1.} {\sl For arbitrary pairs $(r,a),(s,b)$, with $0\le r<s$, we have
$$ P^{s,b}_{r,a}\not=\emptyset\quad\hbox{\rm  if and only if }\quad |b-a|\le c(s-r).$$
If $P^{s,b}_{r,a}\not=\emptyset$, then $P^{s,b}_{r,a}$ is either the line segment connecting $(r,a)$ to $(s,b)$, 
or $P^{s,b}_{r,a}$ is a nondegenerate parallelogram containing this line segment.}
\medskip

\noindent
{\it Proof.} The proof of this proposition is routine. \vrule height 6pt width 5pt depth 4pt
\medskip

\proclaim Definition 2.2. Assume that $P^{s,b}_{r,a}\not=\emptyset$. Let $\displaystyle{u={1\over 2}(r+s)}$. Define $I^{s,b}_{r,a}$ to be the projection of the
following set onto the $x$-axis.
$$\{\,(u,x)\,|\,-\infty < x < +\infty\,\}\cap P^{s,b}_{r,a}.$$
Because $ P^{s,b}_{r,a}$ is a parallelogram containing the line segment joining $(r,a)$ to $(s,b)$, and because $r<u<s$, we see that $I^{s,b}_{r,a}$ is 
either a point or a nondegenerate closed interval.

\medskip

\noindent
{\bf Proposition 2.2.} {\sl The intervals  $I^{s,b}_{r,a}$ are given by}
$$I^{s,b}_{r,a}=\cases{ [\,b-{1\over 2}c(s-r),a+{1\over 2}c(s-r)\,],&if $a\le b$;\cr
                        \mathstrut\cr
                        [\,a-{1\over 2}c(s-r),b+{1\over 2}c(s-r)\,],&if $b\le a$.\cr}$$

\medskip

\noindent
{\it Proof.} The proof of this proposition is clear. \vrule height 6pt width 5pt depth 4pt
\medskip

\proclaim Definition 2.3. We shall assume that for each pair $(r,a)$, $(s,b)$ of points in the plane such that $0\le r<s$ and 
$P^{s,b}_{r,a}\not=\emptyset$, we
are given a continuous function
$$\lambda^{s,b}_{r,a}:[0,1]\to I^{s,b}_{r,a}$$
mapping $[0,1]$ onto $I^{s,b}_{r,a}$.
We shall also assume that the mapping $(a,b,r,s,\xi)\mapsto \lambda^{s,b}_{r,a}(\xi)$ is continuous on the set $D_\lambda$, where 
$$D_\lambda=\{\,(a,b,r,s,\xi)\in\hbox{\b R}^5\,|\,0\le r<s,\quad |b-a|\le c(s-r),\quad  and \quad \xi\in [\,0,1\,]\,\}.$$
\medskip

\noindent
{\bf Proposition 2.3.} {\sl Assume that $P^{s,b}_{r,a}\not=\emptyset$. Let $\displaystyle{u={1\over 2}(r+s)}$. If $d=\lambda^{s,b}_{r,a}(\xi)$, 
where $\xi\in [0,1]$ is arbitrary, then}
$$\emptyset \not= P^{u,d}_{r,a},\;\;P^{s,b}_{u,d} \subseteq P^{s,b}_{r,a}.$$
\medskip

\noindent
{\it Proof.} We will only prove that $\emptyset \not= P^{u,d}_{r,a} \subseteq P^{s,b}_{r,a}$ and $P^{u,d}_{r,a}\subseteq P^{s,b}_{r,a}$. 
The rest of the proof is similar. Also, we only give the proof
for the case where $a\le b$, the proof for the case $b\le a$ is similar. We the have
$$I^{s,b}_{r,a}=[\,b-{1\over 2}c(s-r), a+{1\over 2}c(s-r)\,].$$
Because $d=\lambda^{s,b}_{r,a}(\xi)$, it follows from Definition 2.3 that $d\in I^{s,b}_{r,a}$, therefore
$$b-{1\over 2}c(s-r) \le d \le  a+{1\over 2}c(s-r).\eqno(1)$$
By Proposition 2.1, to prove that $P^{u,d}_{r,a}\not=\emptyset$, it suffices to prove that
$$|d-a|\le c(u-r).\eqno(2)$$
Note that $\displaystyle{u-r={1\over 2}(s-r)}$. Hence, from (1), 
$$d\le  a+{1\over 2}c(s-r)=a+c(u-r).$$ 
It also follows from (1) that 
$$d \ge b-{1\over 2}c(s-r)=b-c(u-r)\ge a-c(u-r).$$
We conclude that
$$a-c(u-r) \le d \le a+c(u-r).$$
Therefore, $|d-a| \le c(u-r)$, hence (2) holds. To prove that $P^{u,d}_{r,a}\subseteq P^{s,b}_{r,a}$, let $(t,x)\in P^{u,d}_{r,a}$ be arbitrary. 
Then by definition, $r\le t\le s$ and
$$\eqalign{a-c(t-r)&\le x\le a+c(t-r),\cr
           d-c(u-t)&\le x\le d+c(u-t).\cr}\eqno(3)$$
We have $u-r-s-u$, and by (2), $d\le a+c(u-r)$, consequently, (3) implies that
$$\eqalign{x&\le d+c(u-t)\cr
            &\le [\,a+c(u-r)\,]+c(u-t)\cr
	    &=   a+c(s-u)+c(u-t)\cr
	    &\le b+c(s-u)+c(u-r)\cr
	    &=b+c(s-r).\cr}$$
Thus, $x\le b+c(s-r)$. Similarly, we see that $b-c(s-r)\le x$, therefore, 
$$b-c(s-r) \le x\le b+c(s-r).$$
Hence, by definition, $(t,x)\in B_{s,b}$. The same type of argument shows that $(t,x)\in F_{r,a}$, consequently, 
$(t,x)\in F_{r,a}\cap B_{s,b}=P^{s,b}_{r,a}$. This proves that $P^{u,d}_{r,a}\subseteq P^{s,b}_{r,a}$. \vrule height 6pt width 5pt depth 4pt
\medskip

\proclaim Definition 2.4. Let $0\le r<s$ be arbitrary, and let $a,b$ be real numbers such that $P^{s,b}_{r,a}\not=\emptyset$. Define
$$t_{nj}=r+{j\over 2^n}(s-r),\;\; 0\le j\le 2^n,\;\; n=0,1,2,\ldots.$$
For $n\ge 1$, define
$$\eqalign{V^s_{r,n}&=\{\,t_{nj}\;|\;0<j<2^n\,\},\cr
           V^s_r&=\bigcup\limits^\infty_{n=1}V^s_{r,n},\cr
	   \Omega^s_r&=[0,1]^{V^s_r}.\cr}$$
\medskip

\proclaim Theorem 2.1. Assume that $P^{s,b}_{r,a}\not=\emptyset$. Then there exists a function 
$$\varphi^{s,b}_{r,a}:\Omega^s_r\to \hbox{\b R}^{V^s_r}$$
such that for each $\omega\in\Omega^s_r$, $\varphi^{s,b}_{r,a}(\omega)$ satisfies conditions (a)--(c) below. For each $n\ge 1$,
we define 
$$a_{mj}(\omega)=\varphi^{s,b}_{r,a}(\omega)(t_{mj}),\;\;0\le j\le 2^m,\;\;0\le m\le n.$$
Write 
$$a_{mj}=a_{mj}(\omega), \;\; \varphi(\omega)=\varphi^{s,b}_{r,a}(\omega).$$
\noindent
The properties that $\varphi(\omega)$ has are as follows.
\item{(a)} For all $1\le m\le n$ and $1\le j\le 2^{m-1}$,
$$\emptyset\not=P^{t_{m,2j-1},a_{m,2j-1}}_{t_{m-1,j-1},a_{m-1,j-1}},\;\;P^{t_{m-1,j},a_{m-1,j}}_{t_{m,2j-1},a_{m,2j-1}}\subseteq
                P^{t_{m-1,j},a_{m-1,j}}_{t_{m-1,j-1},a_{m-1,j-1}}.$$
\item{(b)} For all $1\le m\le n$ and $1\le j\le 2^{m-1}$,
$$\varphi(\omega)(t_{m,2j-1})=\lambda^{t_{m-1,j},a_{m-1,j}}_{t_{m-1,j-1},a_{m-1,j-1}}(\omega_{t_{m,2j-1}}).$$
\item{(c)} $P^{t_{nj},a_{nj}}_{t_{n,j-1},a_{n,j-1}}\not=\emptyset,\quad 1\le j\le 2^n.$
\medskip

\noindent
{\it Proof.} Fix $\omega\in\Omega^s_r$. We will use induction on $n$ to define $\varphi(\omega)=\varphi^{s,b}_{r,a}(\omega)$ 
consistently on each $V_n=V^s_{r,n}$.

Define $\varphi(\omega)$ on $V_1={t_{11}}$ as follows. First, define 
$$\varphi(\omega)(r)=a,\quad \varphi(\omega)(s)=b.$$
Because $P^{s,b}_{r,a}\not=\emptyset$, Definition 2.3 gives a function
$$\lambda^{s,b}_{r,a}:[\,0,1\,]\to I^{s,b}_{r,a}.$$
Define $\varphi(\omega)$ on $V_1$ by
$$\varphi(\omega)(t_{11})=\lambda^{s,b}_{r,a}(\omega_{t_{11}}).$$
Define $\displaystyle{u=t_{11}={1\over 2}(r+s)}$. Then Proposition 2.3 implies that, with $d=\lambda^{s,b}_{r,a}(\omega_{t_{11}})$, 
$$\emptyset\not=P^{u,d}_{r,a},\;\; P^{s,b}_{u,d}\subseteq P^{s,b}_{r,a}.$$
It is now easy to see that (a)--(c) hold for $n=1$.

Now assume that $\varphi(\omega)$ has been defined on $V_n$ in such a way that (a)--(c) hold, where $n\ge 1$ is given. We will then define 
$\varphi(\omega)$ on $V_{n+1}$ in such a way that (a)--(c) hold when $n$ is replaced by $n+1$. That is, we want to define $\varphi(\omega)$ on $V_{n+1}$ 
so that the following conditions hold.
\item{$(a')$} For all $1\le m\le n+1$ and $1\le j\le 2^{m-1}$,
$$\emptyset\not=P^{t_{m,2j-1},a_{m,2j-1}}_{t_{m-1,j-1},a_{m-1,j-1}},\;\;P^{t_{m-1,j},a_{m-1,j}}_{t_{m,2j-1},a_{m,2j-1}}\subseteq
                P^{t_{m-1,j},a_{m-1,j}}_{t_{m-1,j-1},a_{m-1,j-1}}.$$
\item{$(b')$} For all $1\le m\le n+1$ and $1\le j\le 2^{m-1}$,
$$\varphi(\omega)(t_{m,2j-1})=\lambda^{t_{m-1,j},a_{m-1,j}}_{t_{m-1,j-1},a_{m-1,j-1}}(\omega_{t_{m,2j-1}}).$$
\item{$(c')$} $P^{t_{n+1,j},a_{n+1,j}}_{t_{n+1,j-1},a_{n+1,j-1}}\not=\emptyset,\quad 1\le j\le 2^{n+1}$.

To define $\varphi(\omega)$ on $V_{n+1}$, let $0<k<2^{n+1}$ be arbitrary. Then $k$ has one of the following forms. 
$$k=\cases{2j,& $1\le j<2^n$;\cr
           \mathstrut\cr
           2j-1,& $1\le j\le 2^n$.\cr}$$
In  case, $k=2j$, we have $t_{n+1,k}=t_{n+1,2j}=t_{nj}$, and hence we define 
$$a_{n+1,k}=\varphi(\omega)(t_{n+1,k})=\varphi(\omega)(t_{nj}).$$
Suppose that $k=2j-1$. By $(c)$ of the inductive hypothesis, we have
$$P^{t_{nj},a_{nj}}_{t_{n,j-1},a_{n,j-1}}\not=\emptyset.$$
Hence, by Definition 2.3, we are given a function
$$\lambda^{t_{nj},a_{nj}}_{t_{n,j-1},a_{n,j-1}}:[\,0,1\,]\to I^{t_{nj},a_{nj}}_{t_{n,j-1},a_{n,j-1}}.$$
Now define
$$a_{n+1,k}=\varphi(\omega)(t_{n+1,k})=\lambda^{t_{nj},a_{nj}}_{t_{n,j-1},a_{n,j-1}}(\omega_{t_{n+1,k}}).$$
This defines $\varphi(\omega)$ on $V_{n+1}$. 

To prove $(a')$, let $1\le j\le 2^{m-1}$, where $1\le m\le m+1$. Suppose that $m\le n$, then by the inductive hypothesis, $(a)$ holds for $n$. Because
$m\le n$, $(a')$ reduces to $(a)$. On the other hand, suppose that $m=n+1$, and let $1\le j\le 2^{m-1}=2^n$. Define
$$u=t_{n+1,2j-1}={1\over 2}(t_{n,j-1}+t_{nj}).$$ Then Proposition 2.3 implies that, with
$$\eqalign{d&=\lambda^{t_{nj},a_{nj}}_{t_{n,j-1},a_{n,j-1}}(\omega_{t_{n+1,k}})\cr
            &=\varphi(\omega)(t_{n+1,2j-1})\cr
	    &=a_{n+1,2j-1},}$$
we have
$$\emptyset\not=P^{u,d}_{t{n,j-1},a_{n,j-1}},\;\; P^{t_{nj},a_{nj}}_{u,d}\subseteq P^{t_{nj},a_{nj}}_{t{n,j-1},a_{n,j-1}}.$$
This statement is equivalent to $(a')$ for the case where $m=n+1$. We conclude that $(a')$ holds for all $1\le m\le n+1$.

To prove $(b')$, let $1\le m\le n+1$, and let $1\le j\le 2^{m-1}$. If $m\le n$, then by the induction hypothesis, $(b)$ holds, and hence $(b')$ holds 
because $(b')$ reduces to $(b)$. On the other hand, suppose that $m=n+1$. Then by definition,
$$\eqalign{\varphi(\omega)(t_{m,2j-1})&=\varphi(\omega)(t_{n+1,2j-1})\cr
                                      &=\lambda^{t_{nj},a_{nj}}_{t_{n,j-1},a_{n,j-1}}(\omega_{t_{n+1,2j-1}})\cr
				      &=\lambda^{t_{m-1,j},a_{m-1,j}}_{t_{m-1,j-1},a_{m-1,j-1}}(\omega_{t_{m,2j-1}}),\cr}$$
which is $(b')$ for the case $m=n+1$. We conclude that $(b')$ holds for all $1\le m\le n+1$.

To prove $(c')$, let $1\le k\le 2^{n+1}$. Then $k$ has one of the following forms.
$$k=\cases{2j,& $1\le j<2^n$;\cr
           \mathstrut\cr
         2j-1,& $1\le j\le 2^n$.\cr}$$
Suppose that $k=2j$. We proved above that $(a')$ holds, hence, with $m=n+1$ in $(a')$, we have
$$P^{t_{n+1,k},a_{n+1,k}}_{t_{n+1,k-1},a_{n+1,k-1}}=P^{t_{nj},a_{nj}}_{t_{n+1,2j-1},a_{n+1,2j-1}}\not=\emptyset.$$
On the other hand, suppose that $k=2j-1$. Then by $(a')$, with $m=n+1$, we have
$$P^{t_{n+1,k},a_{n+1,k}}_{t_{n+1,k-1},a_{n+1,k-1}}=P^{t_{n+1,2j-1},a_{n+1,2j-1}}_{t_{n,j-1},a_{n,j-1}}\not=\emptyset.$$
Because $1\le k\le 2^{n+1}$ is arbitrary, we conclude that $(c')$ holds. Therefore, $(a')-(c')$ hold, and hence we have completed the inductive
definitions of $\varphi(\omega)$ on each $V_n$ in such a way that $(a)$--$(c)$ hold for each $n$.

To complete the definition of $\varphi(\omega)$ as a function on 
$$V^s_r=\bigcup\limits^\infty_{n=1}V^s_{r,n},$$
we show that $\varphi(\omega)$ is consistently defined on $V^s_r$. To this end, not first that by the above construction, $\varphi(\omega)$ has the 
property that for all $n\ge 1$, $\varphi(\omega)(t_{n+1,2j})=\varphi(\omega)(t_{nj})$, for $0<j<2^n$. Now let $1\le m\le n$, and let $t_{mj}\in V_m$, 
that is, let $0<j<2^m$. Then $t_{mj}=t_{n,j2^{n-m}}\in V_n$. Hence we have
$$\varphi(\omega)(t_{n,j2^{n-m}})=\varphi(\omega)(t_{n-1,j2^{n-1-m}})=\cdots=\varphi(\omega)(t_{mj}).$$
This shows that $\varphi(\omega)$ is a well defined function on $V^s_r$. This completes the proof of the theorem. \vrule height 6pt width 5pt depth 4pt

\medskip
\proclaim Theorem 2.2. Assume that $P^{s,b}_{r,a}\not=\emptyset$, where $0\le r<s$. For any $\omega\in \Omega^s_r$, the function 
$$\varphi^{s,b}_{r,a}(\omega):V^s_r\to\hbox{\b R}^{V^s_r}$$
satisfies 
$$\big|\varphi^{s,b}_{r,a}(\omega)(v)-\varphi^{s,b}_{r,a}(\omega)(u)\big|\le c|v-u|, \eqno (1)$$
for all $u,v\in V^s_r$.

\medskip
\noindent{\it Proof.} Fix $\omega\in\Omega^s_r$, and let $x=\varphi^{s,b}_{r,a}(\omega)=\varphi(\omega)$. Since 
$V^s_r=\bigcup\limits^\infty_{n=1}V^s_{r,n}$, we prove $(1)$ by induction on $n$. For $n=1$, $V^s_{r,1}={\,t_{11}\,}$, and hence $(1)$ holds for 
$n=1$. Now assume that $(1)$ holds for $u,v\in V^s_{r,n}$, where $n\ge 1$. Then we want to prove $(1)$ for $u,v\in V^s_{r,n+1}$. To this end, let
$u=t_{n+1,p}$ and $v=t_{n+1,q}$ be arbitrary members of $V^s_{r,n+1}$, where $0<p<q<2^{n+1}$. Then $p$ has the form
$$p=\cases{2j,& $1\le j<2^n$;\cr       
           \mathstrut\cr
           2j-1,& $1\le j\le 2^n$,\cr}$$
and $q$ has the form
$$q=\cases{2k,& $1\le k<2^n$;\cr
           \mathstrut\cr
           2k-1,& $1\le k\le 2^n$.\cr}$$
Therefore we must consider the following cases.
\item{$(a)$} $p=2j$ and $q=2k$, where $0<j,k<2^n$, and $j\le k$. 
\item{$(b)$} $p=2j$ and $q=2k-1$, where $0<j,k<2^n$, $0<k\le 2^n$, and $2j\le 2k-1$.
\item{$(c)$} $p=2j-1$ and $q=2k$, where $0<j\le 2^n$, $0<k<2^n$, and $2j-1\le 2k$.
\item{$(d)$} $p=2j-1$ and  $q=2k-1$, $0<j\le 2^n$,  $0<k\le 2^n$, and $2j-1\le 2k-1$.

Assume that $(a)$ holds. We have $u=t_{n+1,p}=t_{n+1,2j}=t_{nj}$ and $v=t_{n+1,q}=t_{n+1,2k}=t_{nk}$, and $u,v\in V^s_{r,n}$, therefore $(1)$ holds
by the induction hypothesis. 

Now suppose that $(b)$ holds. Then because $1\le k\le 2^n$, Theorem 2.1(a) implies that 
$$P^{t_{n+1,2k-1},a_{n+1,2k-1}}_{t_{n,k-1},a_{n,k-1}}\not=\emptyset.$$
Hence by Proposition 2.1, we have
$$|x(v)-x(t_{n,k-1})|\le c(v-t_{n,k-1}).\eqno (2)$$
Condition $(b)$ implies that $j\le k-1$. If $j=k-1$, then by $(2)$, we have
$$\eqalign{|x(v)-x(u)|&=|x(v)-x(t_{n,k-1})|\cr
                      &\le c(v-t_{n,k-1})\cr
		                            &=c|v-u|.}$$
On the other hand, suppose that $j<k-1$. By the induction hypothesis, we have
$$|x(t_{nm})-x(t_{n,m-1})|\le c(t_{nm}-t_{n,m-1}),$$
for $j\le m\le k-1$. Therefore, $(2)$ implies that
$$\eqalign{|x(v)-x(u)|&=|[x(v)-x(t_{n,k-1})]+\cdots+[x(t_{n,j+1})-x(t_{nj})]|\cr
                      &\le |x(v)-x(t_{n,k-1})|+\cdots+|x(t_{n,j+1})-x(t_{nj})|\cr
		      &\le c(v-t_{n,k-1})+\cdots+c(t_{n,j+1}-t_{nj})\cr
		      &=c(v-u).\cr}$$
Hence, $(b)$ implies (1).

 Similar arguments show that $(1)$ holds if either $(c)$ or $(d)$ is true. Hence, $(1)$ holds for all $u,v\in V^s_{r,n+1}$. It follows by induction that
$(1)$ is true for all $u,v\in V^s_r$. \vrule height 6pt width 5pt depth 4pt

\medskip
\proclaim Corollary 2.1. Under the hypothesis of Theorem 2.2, for any $\omega\in \Omega^s_r$, we have
$$\big|\varphi^{s,b}_{r,a}(\omega)(v)-\varphi^{s,b}_{r,a}(\omega)(u)\big|\le c|v-u|, \eqno (1)$$
for all $u,v\in V^s_r\cup \{\,r,s\,\}$.

\medskip
\noindent{\it Proof.} Set $\varphi(\omega)=\varphi^{s,b}_{r,a}(\omega)$.Let $u,v\in V^s_r\cup \{\,r,s\,\}$. We consider the following cases.
\item{$(a)$} $u=r$, $v=s$.
\item{$(b)$} $u,v\in V^s_r$.
\item{$(c)$} $u=r$, $v\in V^s_r$.
\item{$(d)$} $u\in V^s_r$, $v=s$.

\noindent
Assume that $(a)$ holds. Because $P^{s,b}_{r,a}\not=\emptyset$, we have $|b-a|\le c(s-r)$, i.e.,
$$|\varphi(\omega)(s)-\varphi(\omega)(r)|\le c|v-u|,$$
which is (1) for $u=r$, $v=s$. If $(b)$ holds, then (1) follows from Theorem 2.2. Suppose that $(c)$ holds. Let $v\in V^s_{r,n}$, say $v=t_{nj}$, 
$0<j<2^n$. Then $u=r=t_{n,0}$. By Theorem 2.1(c), we have $P^{t_{n,1},a_{n,1}}_{t_{n,0},a_{n,0}}\not=\emptyset$, hence Proposition 2.1 implies that
$$|\varphi(\omega)(t_{n,1})-\varphi(\omega)(t_{n,0})|\le c(t_{n,1}-t_{n,0}),$$
that is,
$$|\varphi(\omega)(t_{n,1})-\varphi(\omega)(r)|\le c(t_{n,1}-r).\eqno (2)$$
Then by Theorem 2.2 and (2), we get
$$\eqalign{|\varphi(\omega)(v)-\varphi(\omega)(u)|&=|\varphi(\omega)(t_{nj})-\varphi(\omega)(r)|\cr
                                                 &\le |\varphi(\omega)(t_{nj})-\varphi(\omega)(t_{n,j-1})|+\cdots
						 +|\varphi(\omega)(t_{n,1})-\varphi(\omega)(r)|\cr
						 &\le c(t_{nj}-t_{n,j-1})+\cdots+c(t_{n,1}-r)\cr
						 &=c(t_{nj}-r)\cr
						 &=c|v-u|,\cr}$$
which is (1) for the case where  $u=r$, $v\in V^s_r$. The proof of (1) for the case $(d)$ is similar to the proof of (1) for the case $(c)$.
\vrule height 6pt width 5pt depth 4pt
\medskip

\proclaim Definition 2.5. For $P^{s,b}_{r,a}\not=\emptyset$, define
$$\eqalign{C^s_r&=\{\,x:[r,s]\to\hbox{\b R}\,|\,x\quad\hbox{is continuous on}\quad [\,r,s\,]\,\},\cr
          C^{s,b}_{r,a}&=\{\,x\in C^s_r\,|\,x(r)=a,x(s)=b\,\},\cr
	  \Lambda^s_r&=\{\,x\in C^s_r\,|\,|x(v)-x(u)|\le c|v-u|, \quad\hbox{for all}\quad u,v\in [\,r,s\,]\,\},\cr
     \Lambda^{s,b}_{r,a}&= C^{s,b}_{r,a}\cap \Lambda^s_r,\cr
     L^s_r&=\{\,(a,b)\in\hbox{\b R}^2\,|\,|b-a|\le c(s-r)\,\}.\cr}$$

\medskip
\proclaim Proposition 2.4. Assume that $P^{s,b}_{r,a}\not=\emptyset$. Then for any $\omega\in\Omega^s_r$, the function
$$\varphi^{s,b}_{r,a}(\omega):V^s_r\to\hbox{\b R}^{V^s_r}$$
can be uniquely extended to a function $\varphi^{s,b}_{r,a}(\omega)\in\Lambda^{s,b}_{r,a}$.

\medskip
\noindent{\it Proof.} Fix $\omega\in\Omega$, and set $\varphi^{s,b}_{r,a}(\omega)=\varphi(\omega)$. According to Corollary 2.1, we have
$$|\varphi(\omega)(v)-\varphi(\omega)(u)|\le c|v-u|, \eqno (1)$$
for all $u,v\in V^s_r\cup \{\,r,s\,\}$. To define $\varphi(\omega)$ on all of $[\,r,s\,]$, let $u\in[\,r,s\,]$ be arbitrary. The set 
$V^s_r\cup \{\,r,s\,\}$ is dense in $[\,r,s\,]$, hence there exists a sequence $(u_n)$ in  $[\,r,s\,]$ such that $\lim\limits_{n\to\infty}u_n=u$.
Then (1) implies that the sequence $\big(\varphi(\omega)(u_n)\big)$ is a Cauchy sequence in {\b R}, therefore we may define 
$$\varphi(\omega)(u)=\lim\limits_{n\to\infty}\varphi(\omega)(u_n).$$ 
It is easy to see that (1) then implies that
$$|\varphi(\omega)(v)-\varphi(\omega)(u)|\le c|v-u|,$$
for all $u,v\in[\,r,s\,]$. Thus, the extended function $\varphi(\omega):V^s_r\to\hbox{\b R}^{V^s_r}$ is in $\Lambda^{s,b}_{r,a}$. It is clear from the
definition of $\varphi(\omega)$ on  $[\,r,s\,]$ that this extension is unique. 
\vrule height 6pt width 5pt depth 4pt

\medskip
\proclaim Lemma 2.1. Assume that $P^{s,b}_{r,a}\not=\emptyset$. Let $\displaystyle{u={1\over 2}(r+s)}$. Assume that 
$d\in\hbox{\b R}$ satisfies the following conditions.
$$\eqalignno{|d-a|&\le c(u-r);&(1)\cr
           |d-b|&\le c(s-u).\cr}$$
Then $d\in I^{s,b}_{r,a}$.

\noindent{\it Proof.} We give a proof for the case where $a\le b$. The proof for the case $a\ge b$ is similar. By Proposition 2.2,
$$I^{s,b}_{r,a}=[\,b-{1\over 2}c(s-r),a+{1\over 2}c(s-r)\,].\eqno (2)$$
Hence (1) implies that the following conditions hold.
$$\eqalignno{a-{c\over 2}(s-r)&\le d\le a+{c\over 2}(s-r);&(3)\cr
           b-{c\over 2}(s-r)&\le d\le b+{c\over 2}(s-r).\cr}$$
Then by (2) and (3), we get that $d\in I^{s,b}_{r,a}$.
\vrule height 6pt width 5pt depth 4pt

\medskip
\proclaim Theorem 2.3. For $P^{s,b}_{r,a}\not=\emptyset$, the function
$$\varphi^{s,b}_{r,a}:\Omega^s_r\to\Lambda^{s,b}_{r,a}$$
is onto.

\medskip
\noindent{\it Proof.} Let $x\in\Lambda^{s,b}_{r,a}$. For $n\ge 0$, define $a_{nj}=x(t_{nj})$, where $0\le j\le 2^n$. Because $x\in\Lambda^{s,b}_{r,a}$,
for any $n\ge 0$, and for any $1\le j\le 2^n$, we have
$$\eqalign{|a_{nj}-a_{n,j-1}|&=|x(t_{nj})-x(t_{n,j-1})|\cr
                             &\le c(t_{nj}-t_{n,j-1}).\cr}$$
Therefore, for any $n\ge 0$, and for any $1\le j\le 2^n$, Proposition 2.1 implies that 
$$P^{t_{nj},a_{nj}}_{t_{n,j-1},a_{n,j-1}}\not=\emptyset.\eqno (1)$$
Condition (1) and Definition 2.3 then imply that  for any $n\ge 0$, and for any $1\le j\le 2^n$, we are given the surjective function
$$\lambda^{t_{nj},a_{nj}}_{t_{n,j-1},a_{n,j-1}}:[\,0,1\,]\to I^{t_{nj},a_{nj}}_{t_{n,j-1},a_{n,j-1}}.$$
We will construct by induction an $\omega\in\Omega^s_r$ such that for all $n\ge 1$, the following equation holds:
$$x(t_{n,2j-1})=\lambda^{t_{n-1,j},a_{n-1,j}}_{t_{n-1,j-1},a_{n-1,j-1}}(\omega_{t_{n,2j-1}}),\quad 1\le j\le 2^{n-1}.\eqno (2)$$
To prove (2) for $n=1$, note first that because $P^{s,b}_{r,a}\not=\emptyset$, Definition 2.3 gives the surjective function  
$$\lambda^{s,b}_{r,a}:[\,0,1\,]\to I^{s,b}_{r,a}.$$
Because $x\in\Lambda^{s,b}_{r,a}$ and $\displaystyle{t_{11}={1\over 2}(s+r)}$, we have
$$\eqalignno{|x(t_{11})-a|&=|x(t_{11})-x(r)|\le c(t_{11}-r);&(3)\cr
             |b-x(t_{11})|&=|x(s)-x(t_{11})|\le c(s-t_{11}).\cr}$$
Lemma 2.1 and (3) then imply that $x(t_{11})\in I^{s,b}_{r,a}$. Hence there exists $\omega_{t_{11}}\in [\,0,1\,]$ such that 
$$\eqalignno{x(t_{11})&=\lambda^{s,b}_{r,a}(\omega_{t_{11}}) &(4)\cr
                    &=\lambda^{t_{01},a_{01}}_{t_{00},a_{00}}(\omega_{t_{11}}).\cr}$$
Statement (4) gives that (2) holds for $n=1$. Therefore, as the inductive hypothesis, we assume the $n\ge 1$ is given, and that the numbers
$\omega_{t_{nj}}$, $0\le j\le 2^n$ have been constructed in such a way that (2) holds for $n$. We then want to prove (2) for the case where $n$ is
replaced by $n+1$, that is, we want to construct numbers  $\omega_{t_{n+1,k}}$, $0\le k\le 2^{n+1}$ in $[\,0,1\,]$ such that the following statement holds:
$$x(t_{n+1,2j-1})=\lambda^{t_{nj},a_{nj}}_{t_{n,j-1},a_{n,j-1}}(\omega_{t_{n+1,2j-1}}),\quad 1\le j\le 2^n.\eqno (5)$$
To construct  $\omega_{t_{n+1,k}}$, $0\le k\le 2^{n+1}$, let $1\le k\le 2^{n+1}$ be arbitrary. Then $k$ has one of the following forms.
$$k=\cases{2j,& $1\le j\le 2^n$;\cr
           \mathstrut\cr
	   2j-1,& $1\le j\le 2^n$.\cr}$$
If $k=2j$, $1\le j\le 2^n$, define 
$$\omega_{t_{n+1,k}}=\omega_{t_{nj}}.$$
Assume that $k=2j-1$, $1\le j\le 2^n$. Because $x\in\Lambda^{s,b}_{r,a}$ and $\displaystyle{t_{n+1,2j-1}={1\over 2}(t_{n,j-1}+t_{nj})}$, we have
$$\eqalignno{|x(t_{n+1,2j-1})-a_{n,j-1}|&=|x(t_{n+1,2j-1})-x(t_{n,j-1})|\le c(t_{n+1,2j-1}-t_{n,j-1});&(6)\cr
             |a_{nj}-x(t_{n+1,2j-1})|&=|x(t_{nj})-x(t_{n+1,2j-1})|\le c(t_{nj}-t_{n+1,2j-1}).\cr}$$
It follows from (6) and Lemma 2.1 that 
$$x(t_{n+1,2j-1})\in I^{t_{nj},a_{nj}}_{t_{n,j-1},a_{n,j-1}}.$$
Consequently, there exists $\omega_{t_{n+1,2j-1}}\in [\,0,1\,]$ such that 
$$x(t_{n+1,2j-1})=\lambda^{t_{nj},a_{nj}}_{t_{n,j-1},a_{n,j-1}}(\omega_{t_{n+1,2j-1}}). \eqno (7)$$
This defines $\omega_{t_{n+1,k}}=\omega_{t_{n+1,2j-1}}$. It follows from (7) that the numbers $\omega_{t_{n+1,k}}$, $0\le k\le 2^{n+1}$ satisfy (5). This
completes the inductive construction of $\omega\in\Omega^s_r$ such that (2) holds for all $n\ge 1$. We claim that for $n\ge 1$,
$$x(t_{nj})=\varphi^{s,b}_{r,a}(\omega)(t_{nj}),\quad 1\le j<2^n. \eqno (8)$$
We prove (8) by induction. To this end, write $\varphi(\omega)=\varphi^{s,b}_{r,a}(\omega)$. For $n\ge 0$, define
$$b_{nj}=\varphi(\omega)(t_{nj}),\quad 0\le j\le 2^n.$$
According to Theorem 2.1(b), for $n\ge 1$, we have
$$\varphi(\omega)(t_{n,2j-1})=\lambda^{t_{n-1,j},b_{n-1,j}}_{t_{n-1,j-1},b_{n-1,j-1}}(\omega)(t_{n,2j-1}),\quad 1\le j\le 2^{n-1}. \eqno (9)$$
A simple computation shows that 
$$x(t_{11})=\varphi(\omega)(t_{11})=\lambda^{s,b}_{r,a}(\omega_{t_{11}}),$$
and hence (8) holds for $n=1$. Now assume that (8) holds for $n\ge 1$. We want to then prove that (8) holds when $n$ is replaced by $n+1$, i.e.,
$$x(t_{n+1,j})=\varphi(\omega)(t_{n+1,j}),\quad 1\le j<2^{n+1}. \eqno (10)$$
Let $1\le k<2^{n+1}$. Then $k$ has one of the following forms.
$$k=\cases{2j,& $1\le j<2^n$;\cr
           \mathstrut\cr
	   2j-1,& $1\le j\le 2^n$.\cr}$$
If $k=2j$, $1\le j<2^n$, then by the inductive hypothesis, (8) holds, and hence
$$x(t_{n+1,k})=x(t_{n+1,2j})=x(t_{nj})=\varphi(\omega)(t_{nj})=\varphi(\omega)(t_{n+1,2j})=\varphi(\omega)(t_{n+1,k}). \eqno (11)$$
On the other hand, suppose that $k=2j-1$, $1\le j\le 2^n$. By the induction hypothesis, (8) holds, hence we have
$$\eqalignno{ a_{n,j-1}&=x(t_{n,j-1})=\varphi(\omega)(t_{n,j-1})=b_{n,j-1}, &(12)\cr
            a_{nj}&=x(t_{nj})=\varphi(\omega)(t_{nj})=b_{nj}.\cr}$$
By (2), (9), and (12), we get
$$\eqalignno{x(t_{n+1,k})&=x(t_{n+1,2j-1})=\lambda^{t_{nj},a_{nj}}_{ t_{n,j-1},a_{n,j-1}}(\omega_{t_{n+1,2j-1}}) &(13)\cr
                         &=\lambda^{t_{nj},b_{nj}}_{t_{n,j-1},b_{n,j-1}}(\omega_{t_{n+1,2j-1}})\cr
			 &=\varphi(\omega)(t_{n+1,2j-1})\cr
			 &=\varphi(\omega)(t_{n+1,k}).\cr}$$
Statements (11) and (13) prove (10). Hence, by induction, the claim (8) holds for all $n\ge 1$. Statement (8) is equivalent to the following:
$$x(t)=\varphi(\omega)(t),\quad t\in V^s_r.$$
Because $x,\varphi(\omega)\in \Lambda^{s,b}_{r,a}$, we see that (8) implies
$$x(t)=\varphi(\omega)(t),\quad t\in \{\,r,\,s\}\cup V^s_r.$$
The set $ \{\,r,\,\}\cup V^s_r$ is dense in $[\,r,s\,]$, therefore we have
$$x(t)=\varphi(\omega)(t),\quad t\in [\,r,s\,],$$
i.e., 
$$\varphi(\omega)=x.$$
It follows that the function
$$\varphi^{s,b}_{r,a}:\Omega^s_r\to\Lambda^{s,b}_{r,a}$$
is onto. \vrule height 6pt width 5pt depth 4pt

\medskip
\proclaim Lemma 2.2. Let $U$ be any open subset of $C^s_r$ such that 
$$\Lambda^s_r\cap U\not=\emptyset.$$
Let $x_0\in \Lambda^s_r\cap U$. Then there exists a $\delta>0$ such that if 
$$\tau:\quad r=t_0<t_1<\cdots<t_{n-1}<t_n=s$$
is any partition of $[\,r,s\,]$ with 
$$\max_{1\le j\le n}\bigtriangleup t_j<\delta, \eqno (1)$$
then there exist open intervals $I_{t_j}$, $0\le j\le n$, such that 
$$x_0\in  \Lambda^s_r\cap U_\tau\subseteq U, \eqno (2)$$
where
$$U_\tau=\{\,x\in C^s_r\,|\, x(t_j)\in I_{t_j},\quad 0\le j\le n\,\}.$$

\medskip
\noindent{\it Proof.} Because $U$ is open in $C^s_r$ and $x_0\in U$, there exists an $\epsilon>0$ such that if $x\in C^s_r$ with
$||x-x_0||<\epsilon$, then $x\in U$. Define $\displaystyle{\delta={1\over 6c}\epsilon}$, and assume that 
$$\tau:\quad  r=t_0<t_1<\cdots<t_{n-1}<t_n=s$$
is a partition of $[\,r,s\,]$ satisfying condition (1). Define
$$I_{t_j}=\bigg(\,x_0(t_j)-{1\over 2}\epsilon,x_0(t_j)+{1\over 2}\epsilon\,\bigg),\quad 0\le j\le n.$$
Let $x\in \ \Lambda^s_r\cap U_\tau$. Let $t\in [\,r,s\,]$ be arbitrary, say $t_{j-1}\le t\le t_j$, for some $1\le j\le n$. Because 
$x\in U_\tau$, we have $\displaystyle{|x(t_{j-1})-x_0(t_{j-1})|<{1\over 2}\epsilon}$, and consequently (1) implies that

$$\eqalign{|x(t)-x_0(t)|  &\le |x(t_j)-x(t)|+|x(t_j)-x(t_{j-1})|+|x(t_{j-1})-x_0(t_{j-1})|\cr
                          &\qquad+|x_0(t)-x_0(t_{j-1})|\cr
			  &\le c(t_j-t)+c(t_j-t_{j-1})+{1\over 2}\epsilon+c(t-t_{j-1})\cr
			  &\le 3c\bigtriangleup t_j+{1\over 2}\epsilon\cr
			  &\le 3c\max_{1\le i\le n}\bigtriangleup t_j+{1\over 2}\epsilon\cr
			  &<   3c\delta+{1\over 2}\epsilon=\epsilon.\cr}$$
This shows that $||x-x_0||<\epsilon$, i.e., $x\in U$. Therefore (2) holds. \vrule height 6pt width 5pt depth 4pt

\medskip
\proclaim Definition 2.6. Let $n\ge 1$ be an arbitrary positive integer, and define
$$\Omega^s_{r,n}=[0,1]^{V^s_{r,n}}.$$
Let $\omega\in\Omega^s_r$ be arbitrary. For $1\le j<2^n$, define
$$\omega_{t_{nj}}=(\omega)_{t_{nj}}.$$
Define the function $\pi_n:\Omega\to \Omega^s_{r,n}$ by
$$\pi_n(\omega)=(\omega_{t_{n1}},\ldots,\omega_{t_{n,2^n-1}}).$$
\medskip

\medskip
\proclaim Theorem 2.4. Let $0\le r<s$ be fixed real numbers. For each $n\ge 1$ and for all $0\le j\le 2^n$, there exists continuous functions 
$$\theta_{nj}:L^s_r\times\Omega^s_{r,n}\to \hbox{\b R},$$
such that for any $(a,b,\omega)\in L^s_r\times\Omega^s_r$,
$$\varphi^{s,b}_{r,a}(\omega)(t_{nj})=\theta_{nj}(a,b,\pi_n(\omega)),\quad 0\le j\le 2^n.\eqno (1)$$

\medskip
\noindent{\it Proof.} Write $\Omega^s_r=\Omega$, $L^s_r=L$ and $\Omega^s_{r,n}=\Omega_n$. We will prove (1) by induction on $n\ge 1$. 
To prove (1) for $n=1$, let $(a,b,\omega_{t_{11}})\in L\times\Omega_1$ be arbitrary, and define
$$\theta_{1j}:L\times \Omega_1\to \hbox{\b R},\quad 0\le j\le 2$$
by 
$$\eqalignno{\theta_{10}(a,b,\omega_{t_{11}})&=a, &(2)\cr
             \theta_{11}(a,b,\omega_{t_{11}})&=\lambda^{s,b}_{r,a}(\omega_{t_{11}}),\cr
	     \theta_{12}(a,b,\omega_{t_{11}})&=b.\cr}$$
Then for any $(a,b,\omega)\in L\times\Omega$, Theorem 2.1(b) implies that
$$\varphi^{s,b}_{r,a}(\omega)(t_{11})=\lambda^{s,b}_{r,a}(\omega_{t_{11}})=\theta_{11}(a,b,\omega_{t_{11}})=\theta_{11}(a,b,\pi_1(\omega)). \eqno (3)$$
Also, by definition, we have
$$\eqalignno{\varphi(\omega)(t_{10})&=a=\theta_{10}(a,b,\omega_{t_{11}})=\theta_{10}(a,b,\pi_1(\omega)), &(4)\cr
           \varphi(\omega)(t_{12})&=b=\theta_{12}(a,b,\omega_{t_{11}})=\theta_{12}(a,b,\pi_1(\omega)).\cr}$$
The functions $\theta_{10}$ and $\theta_{12}$ are clearly continuous on $L\times\Omega_1$. It follows from Definition 2.3 that the function
$$(a,b,\xi)\mapsto\lambda^{s,b}_{r,a}(\xi)$$
is a continuous function on $L\times\Omega_1$. Hence (3) and (4) together give (1) for the case $n=1$. 
Now assume that $n\ge 1$ is given and
that there exist continuous functions 
$$\theta_{nj}:L\times\Omega_n\to\hbox{\b R},\quad 0\le j\le 2^n,$$
such that (1) holds. We then want to construct continuous functions
$$\theta_{n+1,k}:L\times\Omega_{n+1}\to\hbox{\b R},\quad 0\le j\le 2^{n+1},$$
such that for all $(a,b,\omega)\in L\times\Omega$, the following statement holds:
$$\varphi^{s,b}_{r,a}(\omega)(t_{n+1,k})=\theta_{n+1,k}(a,b,\pi_{n+1}(\omega)),\quad 0\le k\le 2^{n+1}.\eqno (5)$$
To this end, let $(a,b,\omega_{t_{n+1,1}},\ldots,\omega_{t_{n+1,2^{n+1}-1}})\in L\times\Omega_{n+1}$, 
and let $0\le k\le 2^{n+1}$. Then $k$ has one of the following forms:
$$k=\cases{ 2j,& $0\le j\le 2^n$;\cr
             \mathstrut\cr
            2j-1,& $1\le j\le 2^n$.\cr}$$
If $k=2j$, $0\le j\le 2^n$, then define
$$\theta_{n+1,k}(a,b,\omega_{t_{n+1,1}},\ldots,\omega_{t_{n+1,2^{n+1}i-1}})=\theta_{nj}(a,b,\omega_{t_{n1}},\ldots,\omega_{t_{n,2^n-1}}). \eqno (6)$$
On the other hand, assume that $k=2j-1$, $1\le j\le 2^n$, and let $\omega^{n+1}\in\Omega$ be any member of $\Omega$ such that 
$$\pi_{n+1}(\omega^{n+1})=(\omega_{t_{n+1,1}},\ldots,\omega_{t_{n+1,2^{n+1}-1}}).$$
Then we have
$$\pi_n(\omega^{n+1})=(\omega_{t_{n1}},\ldots,\omega_{t_{n,2^n-1}}).$$
By the induction hypothesis, (1) holds for $n$, hence we have
$$\eqalign{\theta_{nj}(a,b,\omega_{t_{n1}},\ldots,\omega_{t_{n,2^n-1}})&=\varphi^{s,b}_{r,a}(\omega^{n+1})(t_{nj}),\cr
           \theta_{n,j-1}(a,b,\omega_{t_{n1}},\ldots,\omega_{t_{n,2^n-1}})&=\varphi^{s,b}_{r,a}(\omega^{n+1})(t_{n,j-1}).\cr}$$
Therefore, we have
$$\eqalign{|\theta_{nj}(a,b,\omega_{t_{n1}},\ldots,\omega_{t_{n,2^n-1}})-\theta_{n,j-1}(a,b,\omega_{t_{n1}},\ldots,\omega_{t_{n,2^n-1}})|\hfill\cr
           \hfill=|\varphi^{s,b}_{r,a}(\omega^{n+1})(t_{nj})-\varphi^{s,b}_{r,a}(\omega^{n+1})(t_{n,j-1})|\cr
         \le c(t_{nj}-t_{n,j-1}).\quad\quad\quad\quad\quad\quad\quad\quad\quad\quad\cr}$$
Consequently, we have
$$(\,\theta_{n,j-1}(a,b,\omega_{t_{n1}},\ldots,\omega_{t_{n,2^n-1}}),\theta_{nj}(a,b,\omega_{t_{n1}},\ldots,\omega_{t_{n,2^n-1}}),t_{n,j-1},t_{nj},
   \omega_{t_{2j-1}}\,)\in D_\lambda. \eqno (7)$$
Therefore, by Definition 2.3, if we write
$$\eqalignno{b_{n,j-1}&=\theta_{n,j-1}(a,b,\omega_{t_{n1}},\ldots,\omega_{t_{n,2^n-1}}), &(8)\cr
            b_{nj}&=\theta_{nj}(a,b,\omega_{t_{n1}},\ldots,\omega_{t_{n,2^n-1}}),\cr}$$
then (7) implies that we may define
$$\theta_{n+1,k}(a,b,\omega_{t_{n+1,1}},\ldots,\omega_{t_{n+1,2^{n+1}-1}})=
\lambda^{t_{nj},b_{nj}}_{t_{n,j-1},b_{n,j-1}}(\omega_{t_{n+1,2j-1}}). \eqno (9)$$
Definitions (6) and (9) together give the definition of $\theta_{n+1,k}$ on $\Omega_{n+1}$, for $0\le k\le 2^{n+1}$. By Definition 2.3, the functions
$$(p,q,u,v,\xi)\to \lambda^{v,q}_{u,p}(\xi)$$
are continuous on $D_\lambda$, and by the induction hypothesis, the functions $\theta_{nj}$, $\theta_{n,j-1}$ are continuous on $L\times\Omega_n$, hence
(6), (7), and (9) together imply that for $0\le k\le 2^{n+1}$, $\theta_{n+1,k}$ is continuous on $L\times\Omega_{n+1}$. To show that (5) holds,
let $(a,b,\omega)\in\Omega$ and let $0\le k\le 2^{n+1}$. Set
$$\pi_{n+1}(\omega)=(\omega_{t_{n+1,1}},\ldots,\omega_{t_{n+1,2^{n+1}-1}}).$$
Then $k$ has one of the following forms:
$$k=\cases{ 2j,& $0\le j\le 2^n$;\cr
             \mathstrut\cr
            2j-1,& $1\le j\le 2^n$.\cr}$$
If $k=2j$, $0\le j\le 2^n$, then by the induction hypothesis and (9) together imply that
$$\eqalignno{\varphi^{s,b}_{r,a}(\omega)(t_{n+1,k})&=\varphi^{s,b}_{r,a}(\omega)(t_{n+1,2j}) &(10)\cr
                                       &=\varphi^{s,b}_{r,a}(\omega)(t_{nj})\cr
				       &=\theta_{nj}(a,b,\omega_{t_{n1}},\ldots,\omega_{t_{n,2^n-1}})\cr
				       &=\theta_{n+1,2j}(a,b,\pi_{n+1}(\omega))\cr
				       &=\theta_{n+1,k}(a,b,\pi_{n+1}(\omega)).\cr}$$                                                                                
				       
On the other hand, if $k=2j-1$, $1\le j\le 2^n$, then by Theorem 2.1(b) and (9) together imply that
$$\eqalignno{\varphi^{s,b}_{r,a}(\omega)(t_{n+1,k})&=\varphi^{s,b}_{r,a}(\omega)(t_{n+1,2j-1}) &(11)\cr
                                       &=\lambda^{t_{nj},a_{nj}}_{t_{n,j-1},a_{n,j-1}}(\omega)(t_{n+1,2j-1})\cr
                                       &=\lambda^{t_{nj},b_{nj}}_{t_{n,j-1},b_{n,j-1}}(\omega)(t_{n+1,2j-1})\cr
				       &=\theta_{n+1,2j-1}(a,b,\omega_{t_{n+1,1}},\ldots,\omega_{t_{n+1,2^{n+1}-1}})\cr
				       &=\theta_{n+1,k}(a,b,\pi_{n+1}(\omega)).\cr}$$
Statements (10) and (11) together prove (5). Therefore the inductive proof of (1) is complete. \vrule height 6pt width 5pt depth 4pt

\medskip
\proclaim Theorem 2.5. Fix $0\le r<s$. Define the function
$$\phi^s_r:L^s_r\times \Omega^s_r\to\Lambda^s_r$$
by
$$\phi^s_r(a,b,\omega)=\varphi^{s,b}_{r,a}(\omega),\quad (a,b,\omega)\in:L^s_r\times \Omega^s_r.$$
Then $\phi^s_r$ is continuous. In particular, for any $(a,b)\in L^s_r$, the function 
$$\varphi^{s,b}_{r,a}:\Omega^s_r\to\Lambda^{s,b}_{r,a}$$
is continuous.
\medskip

\noindent{\it Proof.} Set $\Omega=\Omega^s_r$ and $\phi^s_r=\phi$ . For $n\ge 1$ and $0\le j\le n$, set $\Omega^s_{r,n}=\Omega_n$.
Let $(a_0,b_0,\omega^0)\in L\times\Omega$ be arbitrary, and define
$$x_0=\phi(a_0,b_0,\omega^0).$$
Let $U$ be any open set in $C^s_r$ such that $x_0\in U$. We want to find an open set $X$ in  $L\times\Omega$ such that
$$(a_0,b_0,\omega_0)\in X \hbox{ and }\phi(X)\subseteq U. \eqno (1).$$
To prove (1) let $\delta>0$ (with respect to $U$) be a positive number given in the hypothesis of Lemma 2.2. Select $n\ge 1$ so large that 
$${1\over 2^n}(s-r)<\delta. \eqno (2)$$
Let $\tau_n$ be the partition of $[\,r,s\,]$ defined by
$$r=t_{n,0},\ldots,t_{n,2^n}=s.$$
By Lemma 2.2, (2) implies that there exists open intervals
$$I_{t_{n,0}},\ldots,I_{t_{n,2^n}}$$
such that 
$$x_0\in \Lambda^s_r\cap U_{\tau_n}\subseteq \Lambda^s_r\cap U, \eqno (3)$$
where
$$ U_{\tau_n}=\{\,x\in C^s_r\,|\,x(t_{nj})\in I_{t_{nj}},\quad 0\le j\le 2^n\,\}.$$
By (3), $x_0\in  U_{\tau_n}$, and hence Theorem 2.4 implies that
$$x_0(t_{nj})=\phi(a_0,b_0,\omega^0)(t_{nj})=\varphi^{s,b}_{r,a}(\omega^0)(t_{nj})=
\theta_{nj}(a,b,\pi_n(\omega^0))\in I_{t_{nj}},\quad 0\le j\le 2^n.\eqno (4)$$
By Theorem 2.4, the function $(a,b,\omega)\mapsto\theta_{nj}(a,b,\pi_n(\omega))$ is continuous on $L\times\Omega$, hence (4) implies that for 
each $0\le j\le 2^n$, there exists an open sets $V_{nj}$ in $L$ and $W_{nj}$ in $\Omega$ such that 
$$(a_0,b_0,\omega^0)\in V_{nj}\times W_{nj}\quad {\rm and}\quad \theta_{nj}(a,b,\pi_n(\omega))\subseteq I_{t_{nj}},
\quad (a,b,\omega)\in V_{nj}\times W_{nj}.\eqno (5)$$
Define
$$X=\bigcap\limits^{2^n}_{j=0}\,V_{nj}\times W_{nj}.$$
Then $X$ is open in $L\times\Omega$ and $(a_0,b_0,\omega^0)\in X$. Let $(a,b,\omega)\in X$ and $0\le j\le 2^n$. Then (5) implies that
$$\phi(a,b,\omega)(t_{nj})=\varphi^{s,b}_{r,s}(\omega)(t_{nj})=\theta_{nj}(a,b,\pi_n(\omega))\in I_{t_{nj}}.$$
Consequently, because $0\le j\le 2^n$ is arbitrary, (3) implies that 
$$\phi(a,b,\omega)=\varphi^{s,b}_{r,s}(\omega)\in \Lambda^{s,b}_{r,a}\cap U_{\tau_n}\subseteq\Lambda^s_r\cap U_{\tau_n}\subseteq U.$$
Because $\omega\in W$ is arbitrary, we see that statement (1) holds. Therefore $\phi$ is continuous on $L\times\Omega$. 
\vrule height 6pt width 5pt depth 4pt

\medskip
\proclaim Definition 2.7. Fix $0\le r<s$. For $n\ge 1$, define
$$\eqalign{\buildrel\wedge\over V\mathstrut^s_{r,n}&=V^s_{r,n}\cup\{s\},\cr
           \buildrel\vee\over V\mathstrut^s_{r,n}&=V^s_{r,n}\cup\{r\},\cr
	   \buildrel\wedge\over V\mathstrut^s_r&=\bigcup\limits^\infty_{n=1}\buildrel\wedge\over V\mathstrut^s_{r,n}=V^s_r\cup \{\,s\,\},\cr
	   \buildrel\vee\over V\mathstrut^s_r&=\bigcup\limits^\infty_{n=1}\buildrel\vee\over V\mathstrut^s_{r,n}=V^s_r\cup \{\,r\,\},\cr
	   \buildrel\wedge\over\Omega\mathstrut^s_r&=[\,0,1\,]^{\buildrel\wedge\over V\mathstrut^s_r},\cr
	   \buildrel\vee\over\Omega\mathstrut^s_r&=[\,0,1\,]^{\buildrel\vee\over V\mathstrut^s_r}.\cr}$$
Now define the functions
$$\eqalign{\buildrel\wedge\over\pi:&\buildrel\wedge\over\Omega\mathstrut^s_r\to\Omega^s_r,\cr
           \buildrel\vee\over\pi:&\buildrel\vee\over\Omega\mathstrut^s_r\to\Omega^s_r\cr}$$
as follows. For $\buildrel\wedge\over\omega\,\in\,\buildrel\wedge\over\Omega\mathstrut^s_r$, 
define $\omega=\buildrel\wedge\over\pi(\buildrel\wedge\over\omega)\in\Omega^s_r$ by
$$\omega_t =\, \buildrel\wedge\over\omega_t,$$
where $t\in V^s_r$. For  $\buildrel\vee\over\omega\,\in\,\buildrel\vee\over\Omega\mathstrut^s_r$, 
define $\omega=\buildrel\vee\over\pi(\buildrel\vee\over\omega)\in\Omega^s_r$ by
$$\omega_t =\, \buildrel\vee\over\omega_t,$$
where $t\in V^s_r$. Finally, for arbitrary $a,b\in\hbox{\b R}$, define
$$\eqalign{\Lambda^s_{r,a}&=\{\,x\in\Lambda^s_r\,|\,x(r)=a\,\},\cr
         \Lambda^{s,b}_a&=\{\,x\in\Lambda^s_r\,|\,x(s)=b\,\}.\cr}$$

\medskip
\proclaim Definition 2.8. Let $0\le r<s$ be arbitrary. For arbitrary $a\in\hbox{\b R}$ define
$$I^s_{r,a}=[\,a-c(s-r),a+c(s-r)\,]. $$
We shall assume that for every $a\in\,${\b R}, we are given a continuous function
$$\lambda^s_{r,a}: [0,1]\to I^s_{r,a}$$
onto $I^s_{r,a}$. Moreover, we will assume that the mapping 
$$(a,r,s,\xi)\to \lambda^s_{r,a}(\xi)$$
is continuous on the set 
$$\{\,(a,r,s,\xi)\,|\,a\in\hbox{\b R},0\le r<s,\xi\in [0,1]\,\}.$$
For $a\in\,${\b R}, define
$$\varphi^s_{r,a}:\buildrel\wedge\over\Omega\mathstrut^s_r\to\Lambda^s_{r,a}$$
by
$$\varphi^s_{r,a}(\buildrel\wedge\over\omega)=\varphi^{s,b}_{r,a}(\buildrel\wedge\over\pi(\buildrel\wedge\over\omega)),$$
where $\buildrel\wedge\over\omega\,\in\,\buildrel\wedge\over\Omega\mathstrut^s_r$, and $b=\lambda^s_{r,a}(\buildrel\wedge\over\omega_s)$. For $b\in\,${\b R},
define 
$$\varphi^{s,b}_r:\buildrel\vee\over\Omega\mathstrut^s_r\to\Lambda^{s,b}_r$$
by
$$\varphi^{s,b}_r(\buildrel\vee\over\omega)=\varphi^{s,b}_{r,a}(\buildrel\vee\over\pi(\buildrel\vee\over\omega)),$$
where $\buildrel\vee\over\omega\,\in\,\buildrel\vee\over\Omega\mathstrut^s_r$, and $a=\lambda^s_{r,b}(\buildrel\vee\over\omega_r)$. For $a\in\,${\b R} and
$\buildrel\wedge\over\omega\,\in\,\buildrel\wedge\over\Omega\mathstrut^s_r$, if  $b=\lambda^s_{r,a}(\buildrel\wedge\over\omega_s)$, then $b\in I^s_{r,a}$, and
hence
$$a-c(s-r)\le b\le a+c(s-r),$$
that is 
$$|b-a|\le c(s-r).$$
It follows from Proposition 2.1 that $P^{s,b}_{r,a}\not=\emptyset$ Hence by Proposition 2.4, we see that 
$$\varphi^s_{r,a}(\buildrel\wedge\over\omega)=\varphi^{s,b}_{r,a}(\buildrel\wedge\over\pi(\buildrel\wedge\over\omega))
  \in\Lambda^{s,b}_{r,a}\subseteq\Lambda^s_{r,a}.$$
A similar argument shows that for $b\in\,${\b R} and $\buildrel\vee\over\omega\,\in\,\buildrel\vee\over\Omega\mathstrut^s_r$, if
$a=\lambda^s_{r,b}(\buildrel\wedge\over\omega_s)$, then
$$\varphi^{s,b}_r(\buildrel\vee\over\omega)=\varphi^{s,b}_{r,a}(\buildrel\vee\over\pi(\buildrel\vee\over\omega))
  \in\Lambda^{s,b}_{r,a}\subseteq\Lambda^{s,b}_r.$$
Consequently, for $a,b\in\,${\b R}, we have
$$\eqalign{\varphi^s_{r,a}:&\buildrel\wedge\over\Omega\mathstrut^s_r\to\Lambda^s_{r,a},\cr
           \varphi^{s,b}_r:&\buildrel\vee\over\Omega\mathstrut^s_r\to\Lambda^{s,b}_r.\cr}$$

\medskip

\proclaim Theorem 2.6. Let $0\le r<s$ be arbitrary. Then for any $a,b\in\hbox{\b R}$, the functions
$$\eqalign{\varphi^s_{r,a}:&\buildrel\wedge\over\Omega\mathstrut^s_r\to\Lambda^s_{r,a},\cr
           \varphi^{s,b}_r:&\buildrel\vee\over\Omega\mathstrut^s_r\to\Lambda^{s,b}_r\cr}$$
are onto.
\medskip
																			   \noindent
{\it Proof.} Fix $a,b\in\hbox{\b R}$. Let $x\in\Lambda^s_{r,a}$, and set $d=x(s)$. Then 
$$|d-a|=|x(s)-x(a)|\le c(s-r),$$
hence we have $|d-a|\le c(s-r)$, i.e.,
$$a-c(s-r)\le d\le a+c(s-r).$$
It follows that $d\in I^s_{r,a}$. By Definition 2.8, the function $\lambda^s_{r,a}:[0,1]\to I^s_{r,a}$ is onto, hence there exists 
$\buildrel\wedge\over\omega_s\;\in[0,1]$ such that $d=\lambda^s_{r,a}(\buildrel\wedge\over\omega_s)$. Because $x\in\Lambda^s_{r,a}$, we have
$x\in\Lambda^{s,d}_{r,a}$, therefore Theorem 2.3 implies that there exists $\omega\in\Omega^s_r$ such that $x=\varphi^{s,d}_{r,a}(\omega)$. Define
$\buildrel\wedge\over\omega\,\in\,\buildrel\wedge\over\Omega\mathstrut^s_r$ by
$$\buildrel\wedge\over\omega_t=\cases{\omega_t,&if $t\in V^s_r$;\cr
                      \buildrel\wedge\over\omega_s,&if $t=s$.\cr}$$            
Then $\buildrel\wedge\over\pi(\buildrel\wedge\over\omega)=\omega$, and hence
$$x=\varphi^{s,d}_{r,a}(\omega)=\varphi^{s,d}_{r,a}(\buildrel\wedge\over\pi(\buildrel\wedge\over\omega))=\varphi^s_{r,a}(\buildrel\wedge\over\omega).$$
Because $x\in\Lambda^s_{r,a}$ is arbitrary, this proves that $\varphi^s_{r,a}$ is onto. A similar proof shows that $\varphi^{s,b}_r$ is onto.
\vrule height 6pt width 5pt depth 4pt
\medskip
                                                                                                                                                           
\proclaim Theorem 2.7. Let $0\le r<s$ be arbitrary. Define
$$\eqalign{\Wedge\phi\mathstrut^s_r:\;&\hbox{\b R}\;\times\Wedge\Omega\mathstrut^s_r\to\Lambda^s_r,\cr
           \Vee\phi\mathstrut^s_r:\;&\hbox{\b R}\;\times\Vee\Omega\mathstrut^s_r\to\Lambda^s_r\cr}$$
by
$$\eqalign{\Wedge\phi\mathstrut^s_r(a,\Wedge\omega)&=\varphi^s_{r,a}(\Wedge\omega),
           \quad (a,\Wedge\omega)\in\hbox{\b R}\;\times\Wedge\Omega\mathstrut^s_r,\cr
	   \Vee\phi\mathstrut^s_r(b,\Vee\omega)&=\varphi^{s,b}_r(\Vee\omega),\quad (b,\Vee\omega)\in\hbox{\b R}\;\times\Vee\Omega\mathstrut^s_r.\cr}$$
Then $\Wedge\phi\mathstrut^s_r$ and $\Vee\phi\mathstrut^s_r$ are continuous. In particular, for $a,b\in\;\hbox{\b R}$ arbitrary, the
functions
$$\eqalign{\varphi^s_{r,a}:\;\Wedge\Omega\mathstrut^s_r\to\Lambda^s_{r,a},\cr
           \varphi^{s,b}_r:\;\Vee\Omega\mathstrut^s_r\to\Lambda^{s,b}_r\cr}$$
are continuous.
\medskip

 \noindent
{\it Proof.} For $(a,\Wedge\omega)\in\hbox{\b R}\;\times\Wedge\Omega\mathstrut^s_r$, and define $b=\lambda^s_{r,a}(\Wedge\omega_s)$. Then 
$$\eqalignno{\Wedge\phi\mathstrut^s_r(a,\Wedge\omega)&=\varphi^s_{r,a}(\Wedge\omega) &(1)\cr
                                                   &=\varphi^{s,b}_{r,a}(\Wedge\pi(\Wedge\omega)\cr
						   &=\phi^s_r(a,b,\Wedge\pi(\Wedge\omega))\cr
						   &=\phi^s_r(a,\lambda^s_{r,a}(\Wedge\omega_s),\Wedge\pi(\Wedge\omega)).\cr}$$
By Theorem 2.5 the function $\phi^s_r:L^s_r\;\times\Wedge\Omega\mathstrut^s_r\to\Lambda^s_r$ is continuous, and by Definition 2.8 the function
$$\hbox{\b R}\times [0,1]\ni (d,\xi)\mapsto \lambda^s_{r,d}(\xi)$$
is continuous. Hence we see that the function
$$\hbox{\b R}\;\times\Wedge\Omega\mathstrut^s_r\;\ni (a,\Wedge\omega)\mapsto \phi^s_r(a,\lambda^s_{r,a}(\Wedge\omega_s),\Wedge\pi(\Wedge\omega))$$
is continuous. Therefore, (1) implies that the function  $\Wedge\phi\mathstrut^s_r$ is continuous. A similar argument shows that the function
$\Vee\phi\mathstrut^s_r$ is continuous. \vrule height 6pt width 5pt depth 4pt
\medskip
                                                                                                                                                           
\proclaim Definition 2.9. Let $r\ge 0$ be arbitrary, and let $m$ be the smallest integer such that $m\ge r$. let $a\in\hbox{\b R}$ be arbitrary. Define
$$\eqalign{C_r&=\{\,x:[\,r,+\infty\,)\to\hbox{\b R}\,|\,x\quad\hbox{\sl is continuous on}\quad [\,r,+\infty\,)\,\},\cr
           \Lambda_r&=\{\,x\in C_r\,|\,|x(v)-x(u)|\le c|v-u|,\quad\hbox{\sl for all}\quad u,v\in [\,r,+\infty\,)\,\},\cr
	   \Lambda_{r,a}&=\{\,x\in\Lambda_r\,|\,x(r)=a\,\},\cr
	   \Wedge V_r&=\Wedge V\mathstrut^m_r\cup\bigcup\limits^\infty_{j=m}\Wedge V\mathstrut^{j+1}_j,\cr
	   \Wedge\Omega_r&=[\,0,1\,]^{\Wedge V_r}.\cr}$$
For $\Wedge\omega\;\in\;\Wedge\Omega_r$, define $\Wedge\pi_{r,m}(\Wedge\omega)$ and $\Wedge\pi_{j,j+1}(\Wedge\omega)$ by
$$\eqalign{\Wedge\pi_{r,m}(\Wedge\omega)_t&=\,\Wedge\omega_t,\;t\in \Wedge V\mathstrut^m_r,\cr
           \Wedge\pi_{j,j+1}(\Wedge\omega)_t&=\,\Wedge\omega_t,\;j\ge m,\;t\in\Wedge V\mathstrut^{j+1}_j.\cr}$$
Let $\Wedge\omega\;\in\;\Wedge\Omega_r$ be arbitrary. Define $\varphi_{r,a}(\Wedge\omega)$ by induction as follows. 
$$\eqalign{\varphi_{r,a}(\Wedge\omega)(t)&=\,\varphi^m_{r,a}(\Wedge\pi_{r,m}(\Wedge\omega))(t),\; t\in [\,r,m\,],\cr
           \varphi_{r,a}(\Wedge\omega)(t)&=\,\varphi^{m+1}_{m,b}(\Wedge\pi_{m,m+1}(\Wedge\omega))(t),\;
	   b=\varphi_{r,a}(\Wedge\omega)(m),\;t\in[\,m,m+1\,].\cr}$$
Now assume that $\varphi_{r,a}(\Wedge\omega)$ has been defined on $[j,j+1]$ as above, where $j\ge m$. Then define
$$\varphi_{r,a}(\Wedge\omega)(t)=\,\varphi^{j+2}_{j+1,d}(\Wedge\pi_{j+1,j+2}(\Wedge\omega))(t),\;d=\varphi_{r,a}(\Wedge\omega)(j+1),\;t\in [\,j+1,j+2\,].$$
Because the ranges of the functions $\varphi^m_{r,a}$, $\varphi^{m+1}_{m,b}$, and $\varphi^{j+2}_{j+1,d}$ ($j\ge m$) are, respectively,
$\Lambda^m_{r,a}$, $\Lambda^{m+1}_{m,b}$, and $\Lambda^{j+2}_{j+1,d}$, we see that
$$\varphi_{r,a}:\;\Wedge\Omega_r\to\Lambda_{r,a}.$$
\medskip

\proclaim Theorem 2.8. Let $r\ge 0$ be arbitrary. Give the set $C_r$ the compact-open topology, and give $\Lambda_r\subseteq C_r$ the induced subspace
topology. Define the function
$$\Wedge\phi_r:\hbox{\b R}\;\times\Wedge\Omega_r\to\Lambda_r$$
by 
$$\Wedge\phi_r(a,\Wedge\omega)=\varphi_{r,a}(\Wedge\omega),\quad (a,\Wedge\omega)\in\hbox{\b R}\;\times\Wedge\Omega_r.$$
Then $\Wedge\phi_r$ is continuous. In particular, for fixed $a\in\;\hbox{\b R}$, the function
$$\varphi_{r,a}:\;\Wedge\Omega_r\to\Lambda_{r,a}$$
is continuous.
\medskip

\noindent
{\it Proof.} Note first that the compact-open topology on $C_r$ coincides with the topology of compact convergence on $C_r$ 
(see {\bf [M]: Theorem 5.1}). Recall that a basis for the topology of compact convergence on $C_r$ consists of all sets of the form
$$B_C(x_0,\epsilon)=\{\,x\in C_r\,|\,\sup\limits_{t\in C}|x(t)-x_0(t)|<\epsilon\,\},$$
where $\epsilon >0$ is arbitrary, and $C$ is and arbitrary compact subset of $[r,\infty)$. Therefore, to show that  $\Wedge\phi_r$ is continuous on
$\hbox{\b R}\;\times\Wedge\Omega_r$, it suffices to show that for arbitrary $(a_0,\Wedge\omega\mathstrut^0)\in\;\hbox{\b R}\;\times\Wedge\Omega_r$,
if $x_0=\;\Wedge\phi_r(a_0,\Wedge\omega\mathstrut^0)$, then for every $\epsilon>0$ and for every compact subset $C$ of $[r,\infty)$, there
exists an open subset $W$ of $\hbox{\b R}\;\times\Wedge\Omega_r$ such that 
$$(a_0,\Wedge\omega\mathstrut^0)\in W \quad \hbox{\rm and} \quad \Wedge\phi_r(W)\subseteq B_C(x_0,\epsilon).\eqno (1)$$
To this end, let $m$ be the smallest positive integer greater than $r$. Let $\epsilon>0$ be arbitrary, and let $C$ be any nonempty
compact subset of $[r,\infty)$. Then there exists a positive integer $k\ge m$ such that $C\subseteq [r,k+1]$. We claim that there exist positive
numbers 
$$0<\epsilon_r,\;\epsilon_m,\;\cdots\epsilon_k,\;\epsilon_{k+1}=\epsilon,$$
and open subsets
$$W_r\subseteq\;\Wedge\Omega\mathstrut^m_r,\quad W_m\subseteq\;\Wedge\Omega\mathstrut^{m+1}_m,\ldots,W_k\subseteq\;\Wedge\Omega\mathstrut^{k+1}_k,$$
having the following properties, where
$$I_r=(x_0(r)-\epsilon_r,x_0(r)+\epsilon_r) \quad\hbox{\rm and}\quad I_j=(x_0(j)-\epsilon_j,x_0(j)+\epsilon_j),\quad m\le j\le k+1.$$
$$\displaylines{\hfill \epsilon_r<\epsilon_j,\quad m\le j\le k+1;\hfill\llap{(2)}\cr
                \hfill \epsilon_j<\epsilon_{j+1},\quad m\le j\le k;\hfill\llap{(3)}\cr
                \hfill \Wedge\pi_{r,m}(\Wedge\omega\mathstrut^0)\;\in W_r\quad
		\hbox{\rm and}\quad \Wedge\phi\mathstrut^m_r(I_r\times W_r)\subseteq B_{[r,m]}(x_0,\epsilon_m);\hfill\llap{(4)}\cr
		\hfill \Wedge\pi_{j,j+1}(\Wedge\omega\mathstrut^0)\;\in W_j\quad\hbox{\rm and}\quad
		\Wedge\phi\mathstrut^{j+1}_j(I_j\times W_j)\subseteq B_{[j,j+1]}(x_0,\epsilon_{j+1}), \quad m\le j\le k.\hfill\llap{(5)}\cr}$$
To prove the claim, we first use ``backward induction'' on $j$ to define $\epsilon_j$ and $W_j$ for $m\le j\le k$, in 
such a way that (3) and (5) hold. By Theorem 2.7, the function 
$$\Wedge\phi\mathstrut^{k+1}_k:\hbox{\b R}\;\times\;\Wedge\Omega\mathstrut^{k+1}_k\to\Lambda^{k+1}_k$$
is continuous. Hence there exists an $0<\epsilon_k<\epsilon_{k+1}=\epsilon$ and an open subset $W_k\subseteq\;\Wedge\Omega\mathstrut^{k+1}_k$ such that
$$\Wedge\pi_{k,k+1}(\Wedge\omega\mathstrut^0)\;\in W_k\quad\hbox{\rm and}\quad
                \Wedge\phi\mathstrut^{k+1}_k(I_k\times W_k)\subseteq B_{[k,k+1]}(x_0,\epsilon_{k+1}).\eqno(6).$$
This defines $\epsilon_k$ and $W_k$. Clearly, (6) implies (5) for the case where $j=k$. Also, (3) clearly holds for $j=k$.

Now assume that $m<j\le k$, and that $\epsilon_j$ and $W_j$ have been defined in such a way that (3) and (5) hold for $j$. By Theorem 2.7,
the function 
$$\Wedge\phi\mathstrut^j_{j-1}:\hbox{\b R}\;\times\;\Wedge\Omega\mathstrut^j_{j-1}\to\Lambda^j_{j-1}$$
is continuous, hence there exists an $0<\epsilon_{j-1}<\epsilon_j$ and an open subset $W_{j-1}$ of $\Wedge\Omega\mathstrut^j_{j-1}$ such that
$$\Wedge\pi_{j,j-1}(\Wedge\omega\mathstrut^0)\;\in W_{j-1}\quad\hbox{\rm and}\quad
                \Wedge\phi\mathstrut^j_{j-1}(I_{j-1}\times W_{j-1})\subseteq B_{[j-1,j]}(x_0,\epsilon_j).\eqno(7).$$
This defines $\epsilon_{j-1}$ and $W_{j-1}$. It is clear that (3) holds when $j$ is replaced by $j-1$. Also, statement (7) implies that 
(5) holds when $j$ is replaced by $j-1$. It follows from backward induction on $j$ that (3) and (5) hold for all $m\le j\le k$. 
We now define $\epsilon_r$ and $W_r$ in such a way that (2) and (4) are true. By Theorem 2.7, the function 
$$\Wedge\phi\mathstrut^m_r:\hbox{\b R}\;\times\;\Wedge\Omega\mathstrut^m_r\to\Lambda^m_r$$
is continuous at $(x_0(r),\Wedge\pi_{r,m}(\Wedge\omega\mathstrut^0))$. Hence there exist $0<\epsilon_r<\epsilon_m$ and an open subset 
$W_r$ of $\Wedge\Omega\mathstrut^m_r$ such that 
$$\Wedge\pi_{r,m}(\Wedge\omega\mathstrut^0)\;\in W_r\quad\hbox{\rm and}\quad
                \Wedge\phi\mathstrut^m_r(I_r\times W_r)\subseteq B_{[r,m]}(x_0,\epsilon_m).\eqno(8).$$
This defines $\epsilon_r$ and $W_r$. Because $\epsilon_r<\epsilon_m$, (3) implies (2). Clearly, (8) implies (4). 

To prove (1), define $W$ as follows.
$$W=I_r\times(W_r\times\cdots\times W_k\times\;\Wedge\Omega_{k+1}).$$
Then $W$ is open in $\hbox{\b R}\times\;\Wedge\Omega_r$. By (4) and (5) we see that 
$$\Wedge\pi_{r,m}(\Wedge\omega\mathstrut^0)\in W_r\quad\hbox{\rm and}\quad\Wedge\pi_{j,j+1}(\Wedge\omega\mathstrut^0)\in W_j,\quad m\le j\le k.$$
It follows that 
$$(a_0,\Wedge\omega\mathstrut^0)\in W. \eqno(9)$$
Now let $x\in\;\Wedge\phi_r(W)$ be arbitrary, say 
$$x=\Wedge\phi_r(a,\Wedge\omega),\quad (a,\Wedge\omega)\in W.$$
By definition of $\Wedge\phi_r$, we have
$$x|[r,m]=\Wedge\phi_r(a,\Wedge\omega)|[r,m]=\varphi^m_{r,a}(\Wedge\pi_{r,m}(\Wedge\omega)).$$
Because $\Wedge\pi_{r,m}(\Wedge\omega)\in W_r$ and $x(r)=a\in W_r$, (3) and (4) together imply that 
$$\eqalign{x|[r,m]&=\varphi^m_{r,a}(\Wedge\pi_{r,m}(\Wedge\omega))\cr
                  &=\;\Wedge\phi\mathstrut^m_r(a,\Wedge\pi_{r,m}(\Wedge\omega))\cr
		  &\in\;\Wedge\phi\mathstrut^m_r(I_r\times W_r)\cr
		  &\subseteq B_{[r,m]}(x_0,\epsilon_m)\cr
		  &\subseteq B_{[r,m]}(x_0,\epsilon).\cr}$$
Therefore, we have
$$x|[r,m]\in B_{[r,m]}(x_0,\epsilon_m)\subseteq B_{[r,m]}(x_0,\epsilon). \eqno(10)$$
We claim that 
$$x|[j,j+1]\in B_{[j,j+1]}(x_0,\epsilon_{j+1})\subseteq B_{[j,j+1]}(x_0,\epsilon),\quad m\le j\le k. \eqno(11)$$
We prove this claim by using induction on $m\le j\le k$. By (10), we have 
$$|x(m)-x_0(m)|<\epsilon_m,$$
hence $x(m)\in I_m$. Therefore, by (9),
$$(x(m),\Wedge\pi_{m,m+1}(\Wedge\omega))\;\in I_m\times W_m.$$           
It follows from (5) and Definition 2.9 that, with $b=\varphi_{r,a}(\Wedge\omega)(m)=x(m)$, we have
$$\eqalign{x|[m,m+1]&=\varphi^{m+1}_{m,b}(\Wedge\pi_{m,m+1}(\Wedge\omega))\cr
                    &=\;\Wedge\phi\mathstrut^{m+1}_m(b,\Wedge\pi_{m,m+1}(\Wedge\omega))\cr
		    &\subseteq\;\Wedge\phi\mathstrut^{m+1}_m(I_m\times W_m)\cr
		    &\subseteq B_{[m,m+1]}(x_0,\epsilon_{m+1})\cr
		    &\subseteq B_{[m,m+1]}(x_0,\epsilon).\cr}$$
Therefore (11) holds for $j=m$. Assume that (11) holds for $m\le j<k$. Then by (11), $|x(j+1)-x_0(j+1)|<\epsilon_{j+1}$, and hence $x(j+1)\in I_{j+1}$.
Consequently, by (9), we have
$$(x(j+1),\Wedge\pi_{j+1,j+2}(\Wedge\omega))\in I_{j+1}\times W_{j+1}.$$
It follows from (5) and Definition 2.9 that, with $d=\varphi_{r,a}(\Wedge\omega)(j+1)=x(j+1)$, we have
$$\eqalign{x|[j+1,j+2]&=\varphi^{j+2}_{j+1,d}(\Wedge\pi_{j+1,j+2}(\Wedge\omega))\cr
                    &=\;\Wedge\phi\mathstrut^{j+2}_{j+1}(d,\Wedge\pi_{j+1,j+2}(\Wedge\omega))\cr
		    &\subseteq\;\Wedge\phi\mathstrut^{j+2}_{j+1}(I_{j+1}\times W_{j+1})\cr
                    &\subseteq B_{[j+1,j+2]}(x_0,\epsilon_{j+2})\cr
                    &\subseteq B_{[j+1,j+2]}(x_0,\epsilon).\cr}$$
That is, 
$$x|[j+1,j+2]\in B_{[j+1,j+2]}(x_0,\epsilon_{j+2})\subseteq  B_{[j+1,j+2]}(x_0,\epsilon). \eqno(12)$$
Statement (12) is just statement (11) with $j$ replaced by $j+1$. Hence, by induction, (11) holds. Now, (10) and (11) together imply the following statements.
$$\displaylines{\hfill x|([r,m]\cap C)\in B_{[r,m}(x_0,\epsilon);\hfill\llap{(13)}\cr
              \hfill x|([j,j+1]\cap C)\in B_{[j,j+1]}(x_0,\epsilon),\quad m\le j\le k.\hfill\llap{(14)}\cr}$$
Because $C\subseteq [r,m]\cup\bigcup\limits^k_{j=m}\;[j,j+1]=[r,k+1]$, we see that (13) and (14) together imply that 
$$x\in B_C(x_0,\epsilon).$$
Because $x\in W$ is arbitrary, we see that (1) holds. This completes the proof that 
$$\Wedge\phi_r:\hbox{\b R}\;\times\Wedge\Omega_r\to\Lambda_r$$
is continuous in the compact-open topology on $\Lambda_r$. \vrule height 6pt width 5pt depth 4pt
\medskip
proclaim Theorem 2.9. Let $r\ge 0$ be arbitrary. Then the function
$$\Wedge\phi_r:\hbox{\b R}\;\times\;\Wedge\Omega_r\to\Lambda_r$$
is onto.
\medskip

\noindent
{\it Proof.} Let $m$ be defined as in Theorem 2.8. Let $x\in\Lambda_r$. Then $x|[r,m]\in\Lambda^m_{r,x(r)}$, hence by Theorem 2.6, there exists
$\Wedge\omega\mathstrut^r\in\;\Wedge\Omega\mathstrut^m_r$ such that, with $a=x(r)$, we have
$$x|[r,m]=\varphi^m_{r,a}(\Wedge\omega\mathstrut^r).$$
Now let $j\ge m$ be an arbitrary integer. By Theorem 2.6, there exists an $\Wedge\omega\mathstrut^j\in\;\Wedge\Omega\mathstrut^{j+1}_j$ such that
$$x|[j,j+1]=\varphi^{j+1}_{j,x(j)}(\Wedge\omega\mathstrut^j).$$
Define $\Wedge\omega\in\;\Wedge\Omega_r$ by
$$\Wedge\omega\;=(\Wedge\omega\mathstrut^r,\Wedge\omega\mathstrut^m,\ldots,\Wedge\omega\mathstrut^j,\ldots).$$
We claim that
$$x=\;\Wedge\phi_r(\Wedge\omega).\eqno(1)$$
To prove this claim, note first that by Definition 2.9, for $t\in[r,m]$, we have
$$\eqalignno{x(t)=(x|[r,m])(t)&=\varphi^m_{r,a}(\Wedge\omega\mathstrut^r)(t) &(2)\cr
			    &=\varphi^m_{r,m}(\Wedge\pi_{r,m}(\Wedge\omega))(t)\cr
			    &=\varphi_{r,a}(\Wedge\omega)(t)\cr
			    &=\;\Wedge\phi_r(a,\Wedge\omega)(t).\cr}$$
Define $d=\varphi_{r,x(r)}(\Wedge\omega)(m)$. Then by (2), we have $d=x(r)$. Hence by Definition 2.9, for $t\in[m,m+1]$, we have
$$\eqalignno{x(t)=(x|[m,m+1])(t)&=\varphi^{m+1}_{m,x(m)}(\Wedge\omega\mathstrut^m)(t) &(3)\cr
                              &=\varphi^{m+1}_{m,x(m)}(\Wedge\pi_{m,m+1}(\Wedge\omega))(t)\cr
                              &=\varphi_{r,d}(\Wedge\omega)(t)\cr
                              &=\;\Wedge\phi_r(a,\Wedge\omega)(t).\cr}$$
To finish the proof of the claim (1), we will prove by induction on $j\ge m$ that
$$x(t)=\;\Wedge\phi_r(a,\Wedge\omega)(t),\quad t\in[j,j+1]. \eqno(4).$$
Statement (3) implies that (4) holds for $j=m$. Assume that (4) holds for some $j\ge m$. Define $h=\varphi_{r,a}(\Wedge\omega)(j+1)$. 
Then, by the induction hypothesis, we have 
$$\eqalignno{x(t)=(x|[j+1,j+2])(t)&=\varphi^{j+2}_{j+1,x(j+1)}(\Wedge\omega\mathstrut^{j+1})(t) &(5)\cr
                                &=\varphi^{j+2}_{j+1,h}(\Wedge\pi_{j+1,j+2}(\Wedge\omega))(t)\cr
			        &=\varphi_{r,a}(\Wedge\omega)(t)\cr
			        &=\;\Wedge\phi_r(a,\Wedge\omega)(t).\cr}$$
Statement (5) is obtained from statement (4) by replacing $j$ in (4) by $j+1$. Hence (4) holds by induction on $j\ge m$. Now statements 
(2) and (4) together imply (1). Because $x\in\Lambda_r$ is arbitrary, we have proved that $\Wedge\phi_r$ is onto. \vrule height 6pt width 5pt depth 4pt
\medskip

\proclaim Definition 2.10. Let $0\le r<s$ be arbitrary. Define
$$\eqalign{\Sim\Omega\mathstrut^s_r&=\hbox{\b R}\;\times\;\Wedge\Omega\mathstrut^s_r,\cr
           \Sim\Omega_r&=\hbox{\b R}\;\times\;\Wedge\Omega_r.\cr}$$
Assume that we are given a function 
$$\lambda_r:\hbox{\b R}\to\hbox{\b R}$$
such that $\lambda_r$ is continuous and onto. Then define the functions
$$\eqalign{\varphi^s_r&:\;\Sim\Omega\mathstrut^s_r\to\Lambda^s_r,\cr
           \varphi_r&:\;\Sim\Omega\mathstrut_r\to\Lambda_r\cr}$$
as follows.
$$\eqalign{\varphi^s_r(a,\Wedge\omega)&=\varphi^s_{r,b}(\Wedge\omega),\quad b=\lambda_r(a),\quad (a,\Wedge\omega)\in\;\Sim\Omega\mathstrut^s_r,\cr
           \varphi_r(a,\Wedge\omega)&=\varphi_{r,b}(\Wedge\omega),\quad b=\lambda_r(a),\quad (a,\Wedge\omega)\in\;\Sim\Omega\mathstrut_r.\cr}$$
\medskip

\proclaim Theorem 2.10. Let $0\le r<s$ be arbitrary. Then the functions
$$\eqalign{\varphi^s_r&:\;\Sim\Omega\mathstrut^s_r\to\Lambda^s_r,\cr
           \varphi_r&:\;\Sim\Omega\mathstrut_r\to\Lambda_r\cr}$$
are continuous and onto. 
\medskip

\noindent
{\it Proof.} To prove that $\varphi^s_r$ is onto, let $x\in\Lambda^s_r$. Because $\lambda_r$ is onto, there exists an 
$a\in\hbox{\b R}$ such that $\lambda_r(a)=x(r)$. Define $b=x(r)$, then by Theorem 2.6, there exists $\Wedge\omega\;\in\;\Wedge\Omega\mathstrut^s_r$ such that
$$\eqalign{x&=\varphi^s_{r,b}(\Wedge\omega)\cr
            &=\varphi^s_r(a,\Wedge\omega).\cr}$$
Because $x$ is arbitrary, this shows that $\varphi^s_r$ is onto. A similar argument shows that $\varphi_r$ is onto.

To show that $\varphi^s_r$ is continuous, let $(a,\Wedge\omega)\in\;\Sim\Omega\mathstrut^s_r$. By assumption, $\lambda_r$ is continuous, 
and by Theorem 2.7,
the function
$$\Wedge\phi\mathstrut^s_r:\;\hbox{\b R}\;\times\;\Wedge\Omega\mathstrut^s_r\to\Lambda^s_r$$
is continuous, therefore, for $(s,\Wedge\omega)\in\;\Sim\Omega\mathstrut^s_r$, the function 
$$\eqalign{(a,\Wedge\omega)\mapsto\;\Wedge\phi\mathstrut^s_r(\lambda_r(a),\Wedge\omega)&=\;\Wedge\phi\mathstrut^s_r(b,\Wedge\omega)\cr
                                                                                     &=\varphi^s_{r,b}(\Wedge\omega)\cr
										     &=\varphi^s_r(a,\Wedge\omega).\cr}$$
is continuous. Thus, $\varphi^s_r$ is continuous. A similar argument shows that $\varphi_r$ is continuous. \vrule height 6pt width 5pt depth 4pt

\vskip 1cm
\centerline{\bf 3. CONSTRUCTION OF THE MEASURES $\mu^{\,\square}_{\,\square}$}
\vskip 1cm
\noindent
In this section we use standard results from topology and measure theory to construct the families of regular
Borel measures mentioned in (2) and (3) of the introduction.

\medskip
\proclaim Theorem 3.1 (Ascoli's theorem). Let $X$ be a locally compact Hausdorff space; let $(Y,d)$ be a metric space. 
Let ${\cal C}(X,Y)$ be the space of all
continuous functions from $X$ to $Y$, and consider ${\cal C}(X,Y)$ in the compact-open topology. 
A subset $\cal F$ of ${\cal C}(X,Y)$ has compact closure if and only
if $\cal F$ is equicontinuous and the subset
$${\cal F}_x=\{\,f(x)\,|\,f\in{\cal F}\,\}$$
of $Y$ has compact closure for each $x\in X$.
\medskip

\noindent
{\it Proof.} (See {\bf [M]: Theorem 6.1}.) \vrule height 6pt width 5pt depth 4pt
\medskip

\proclaim Theorem 3.2. Let $0\le r<s$ be arbitrary, and let $(a,b)\in L^s_r$. Then each of the following function spaces is equicontinuous:
$$\Lambda^{s,b}_{r,a},\quad \Lambda^s_{r,a},\quad \Lambda^{s,b}_r,\quad \Lambda_{r,a},\quad \Lambda^s_r,\quad \Lambda_r.$$
\medskip

\noindent
{\it Proof.} We show that $\Lambda^{s,b}_{r,a}$ is equicontinuous. The proof that remaining spaces are equicontinuous is similar. We show that 
$\Lambda^{s,b}_{r,a}$ is equicontinuous at each point $t_0\in[r,s]$. Let $\epsilon>0$ be arbitrary, and set $\displaystyle{\delta={\epsilon\over c}}$.
Then for all $x\in\Lambda^{s,b}_{r,a}$ and all $t\in[r,s]\cap(t_0-\delta,t_0+\delta)$, we have
$$|x(t)-x(t_0)|\le c|t-t_0|<c\delta=\epsilon.$$
Thus, $\Lambda^{s,b}_{r,a}$ is equicontinuous at the arbitrary point $t_0\in [r,s]$, hence $\Lambda^{s,b}_{r,a}$ is equicontinuous. 
\vrule height 6pt width 5pt depth 4pt
\medskip

\proclaim Theorem 3.3. Let $0\le r<s$ be arbitrary. Then for arbitrary $(a,b)\in L^s_r$, $\Lambda^{s,b}_{r,a}$, $\Lambda^s_{r,a}$, $\Lambda^{s,b}_r$,  and 
$\Lambda^s_r$ are compact in the uniform topology. If $a\in\hbox{\b R}$ is arbitrary, then $\Lambda_{r,a}$ is compact in the compact-open topology.
\medskip

\noindent
{\it Proof.} To prove that $\Lambda^{s,b}_{r,a}$ is compact in the the uniform topology, note first that the compact-open topology on 
$\Lambda^{s,b}_{r,a}$ coincides with the topology of compact convergence, and because $[r,s]$ is compact, 
the topology of compact convergence on $\Lambda^{s,b}_{r,a}$
coincides with the uniform topology (see {\bf [M]: Theorem 4.6}). Therefore, it suffices to prove that 
$\Lambda^{s,b}_{r,a}$ is compact in the compact-open topology. To this end, let $t\in [r,s]$, and let $x\in\Lambda^{s,b}_{r,a}$. Then we have
$$\eqalignno{|x(t)|&=|(x(t)-x(r))+x(r)|&(1)\cr
                   &\le |x(t)-x(r)|+|a|\cr
		   &\le c(t-r)+|a|.\cr}$$
Now define ${\cal F}=\Lambda^{s,b}_{r,a}$, and define ${\cal F}_t=\{\,x(t)\,|\,x\in{\cal F}\,\}$. Then (1) implies that ${\cal F}_t$ is bounded, and 
hence it has compact closure. By Theorem 3.2, $\cal F$ is equicontinuous, and hence Ascoli's theorem implies that $\cal F$ has compact closure in the
compact-open topology, that is, ${\cal F}=\Lambda^{s,b}_{r,a}$ is compact in the compact-open topology. A similar argument shows that the spaces 
$\Lambda^s_{r,a}$ and $\Lambda^s_r$ are compact in the uniform topology. 

Now let $a\in\hbox{\b R}$ be arbitrary. Define ${\cal G}=\Lambda_{r,a}$. Let $t\in[r,+\infty)$ be arbitrary, and define 
${\cal G}_t=\{\,x(t)\,|\,x\in{\cal G}\,\}$. Then (1) implies that ${\cal G}_t$ has compact closure. By Theorem 3.2, $\cal G$ is equicontinuous, and
hence by Ascoli's theorem, ${\cal G}=\Lambda_{r,a}$ is compact in the compact-open topology. \vrule height 6pt width 5pt depth 4pt
\medskip

\proclaim Theorem 3.4. Let $0\le r<s$ be arbitrary. Then the space $\Lambda^s_r$ is locally compact in the uniform topology, and the space
$\Lambda_r$ is locally compact in the compact-open topology.
\medskip

\noindent
{\it Proof.} As in the proof of Theorem 3.3, the uniform topology on $\Lambda^s_r$ coincides with the compact-open topology, hence it suffices to
show that $\Lambda^s_r$ is locally compact in the compact-open topology. To this end, let $x_0\in \Lambda^s_r$ be arbitrary. Let $\epsilon >0$ be 
arbitrary, and define $C=\{r\}$, $U=(x_0(r)-\epsilon,x_0(r)+\epsilon)$. Then the set
$$S(C,U)=\{\,x\in \Lambda^s_r\,|\,x(C)\subseteq U\,\}=\{\,x\in \Lambda^s_r\,|\,|x(r)-x_0(r)|<\epsilon\,\}$$
is a basis element in the compact-open topology on  $\Lambda^s_r$, and this basis element contains $x_0$. We claim that $S(C,U)$ has compact closure;
and hence, because $x_0\in \Lambda^s_r$ is arbitrary, $\Lambda^s_r$ is locally compact in the compact-open topology. To prove the claim, let
$x\in S(C,U)$, then we have
$$\eqalignno{|x(r)|&=|(x(r)-x_0(r))+x_0(r)|&(1)\cr
                   &\le |x(r)-x_0(r)|+|x_0(r)|\cr
		   &<\epsilon+|x_0(r)|.\cr}$$
Hence, for any $t\ge r$, we have
$$\eqalignno{|x(t)|&=|(x(t)-x(r))+x(r)|&(2)\cr
                   &\le |x(t)-x(r)|+|x(r)|\cr
		   &< c(t-r)+|x_0(r)|+\epsilon.\cr}$$
Define ${\cal F}=S(C,U)$. Then (2) implies that ${\cal F}_t$ has compact closure. By Theorem 3.2, $\cal F$ is equicontinuous, and hence, because
$t\ge r$ is arbitrary, Ascoli's
theorem gives that ${\cal F}=S(C,U)$ has compact closure in the compact-open topology. Thus,  $\Lambda^s_r$ is locally compact in the 
compact-open topology. A similar argument shows that $\Lambda_r$ is locally compact in the compact-open topology.  \vrule height 6pt width 5pt depth 4pt
\medskip

\proclaim Theorem 3.5. Let $0\le r<s$ be arbitrary. Give the space $\Lambda^s_r$ the uniform topology, and give the space $\Lambda_r$ the compact-open
topology. Assume that the function $\lambda_r:\hbox{\b R}\to\hbox{\b R}$ in Definition 2.10 has the property that if $C$ is any compact subset of 
$\hbox{\b R}$, then $\lambda^{-1}_r(C)$ is compact. Then for any compact subset $A$ of  $\Lambda^s_r$ and any compact subset $B$ of $\Lambda_r$,
the following sets are compact:
$$(\varphi^s_r)^{-1}(A),\quad \varphi^{-1}_r(B).$$
\medskip

\noindent
{\it Proof.} Let $A\subseteq\Lambda^s_r$ be compact. By Theorem 3.2, $\Lambda^s_r$ is equicontinuous, and hence $A$ is also equicontinuous. Ascoli's 
theorem then implies that the set
$$A_r=\{\,x(r)\,|\,x\in A\,\}$$
has compact compact closure. Therefore there exists a $d>0$ such that $|x(r)|\le d$ for all $x\in A$. 
Now let $(a,\Wedge\omega)\in(\varphi^s_r)^{-1}(A)$. Then $\varphi^s_r(a,\Wedge\omega)\in A$. By Definition 2.10, with $b=\lambda_r(a)$, we have
$$\varphi^s_r(a,\Wedge\omega)(r)=\varphi^s_{r,b}(\Wedge\omega)(r)=b=\lambda_r(a).$$
It follows that $|\lambda_r(a)|=|\varphi^s_r(a,\Wedge\omega)(r)|\le d$, i.e., $a\in\lambda^{-1}_r([-d,d])$. The set $\lambda^{-1}_r([-d,d])$ is, by
assumption, compact, hence there exists an $h>0$, depending only on $d$, such that $\lambda^{-1}_r([-d,d])\subseteq [-h,h]$, 
which implies that $a\in [-h,h]$. Because $(a,\Wedge\omega)\in(\varphi^s_r)^{-1}(A)$ is arbitrary, we see that
$$(\varphi^s_r)^{-1}(A)\subseteq [-h,h]\,\times\Wedge\Omega\mathstrut^s_r.$$
By Theorem 2.10, the function $\varphi^s_r$ is continuous, there $(\varphi^s_r)^{-1}(A)$ is a closed subset of the compact set
$[-h,h]\,\times\Wedge\Omega\mathstrut^s_r$. It follows that $(\varphi^s_r)^{-1}(A)$ is compact. Similarly,
$\varphi^{-1}_r(B)$ is compact whenever $B\subseteq\Lambda_r$ is compact in the compact-open topology. \vrule height 6pt width 5pt depth 4pt
\medskip

\proclaim Definition 3.1. Let $X$ and $Y$ be locally compact Hausdorff spaces, and let $\psi$ be a continuous function from $X$ onto $Y$. Let
$(X,{\cal M}_\mu,\mu)$ be a measure space constructed from a nonnegative linear functional on $X$ via the Daniell approach 
(see {\bf [HS]: \S 9}). Assume that one of the following conditions holds.
\item{(a)} $\psi^{-1}(F)$ is compact in $X$ for every compact subset $F$ of $Y$.
\item{(b)} $\mu(X)<+\infty$.
\medskip

\noindent Let ${\cal C}(Y)$ be the space of all continuous complex-valued functions on $Y$, and let ${\cal C}_{00}(Y)$ be the set of all
$f\in{\cal C}(Y)$ such that $f$ has compact support. Then for any $f\in{\cal C}_{00}(Y)$, $f\circ\psi\in {\cal L}_1(X,{\cal M}_\mu,\mu)$, and hence
the mapping
$$J(f)=\int\limits_{X}\,f\circ\psi(x)\,d\mu(x)$$
is a nonnegative linear functional on ${\cal C}_{00}(Y)$. Hence we may construct a measure space $(Y,{\cal M}_\nu,\nu)$ from $J$ via the Daniell 
approach. We then have
$$\int\limits_{Y}\,f(y)\,d\nu(y)=\int\limits_{X}\,f\circ\psi(x)\,d\mu(x),\quad f\in{\cal C}_{00}(Y).$$
The measure $\nu$ is said to be the image of the measure $\mu$ under the continuous function $\psi$.
\medskip

\proclaim Theorem 3.6. The measure $\nu$ constructed in Definition 3.1 has the following properties. 
\item{(a)} For all $\sigma$-finite $\nu$-measurable subsets $B$ of $Y$, 
$$\nu(B)=\mu(\psi^{-1}(B))=\int\limits_{X}\,\chi\circ\psi(x)\,d\mu(x)$$.
\item{(b)} For every $f\in{\cal L}_1(Y,{\cal M}_\nu,\nu)$, $f\circ\psi\in{\cal L}_1(X,{\cal M}_\mu,\mu)$ and
$$\int\limits_{Y}\,f(y)\,d\nu(y)=\int\limits_{X}\,f\circ\psi(x)\,d\mu(x).$$
\medskip

\noindent
{\it Proof.} See {\bf [HS]: Theorem 12.46}. \vrule height 6pt width 5pt depth 4pt
\medskip

\proclaim Theorem 3.7. Let $X$ be a locally compact Hausdorff space and let $\nu$ be a regular measure defined on a $\sigma$-algebra $\cal A$ of
subsets of $X$ such that $(X,{\cal A},\nu)$ is a complete measure space. Suppose that $E\in{\cal A}$ if and only if $E\cap F\in{\cal A}$ for every
compact set $F\subseteq X$. Define $I$ on ${\cal C}_{00}(X)$ by
$$I(f)=\int\limits_{X}\,f(x)\,d\nu(x).$$
Let $(X,{\cal M}_\iota,\iota)$ be the measure space constructed from $I$ via the Daniell approach. Then ${\cal A}={\cal M}_\iota$ and 
$\nu(E)=\iota(E)$ for all $E\in{\cal M}_\iota$.
\medskip

\noindent
{\it Proof.} (See {\bf [HS]: Theorem 12.42}.) \vrule height 6pt width 5pt depth 4pt
\medskip

\proclaim Definition 3.2. Let $0\le r<s$ be arbitrary, and let $(a,b)\in L^s_r$ be arbitrary. Define
$$m^s_r,\quad \Wedge m\mathstrut^s_r,\quad \Vee m\mathstrut^s_r,\quad \Wedge m_r, \eqno (1)$$
to be normalized Lebesgue measure, respectively, on the following product spaces:
$$\Omega^s_r,\quad \Wedge\Omega\mathstrut^s_r,\quad \Vee\Omega\mathstrut^s_r,\quad \Wedge\Omega_r. \eqno (2)$$
Let 
$$\Sim m\mathstrut^s_r,\quad \Sim m_r$$
be Lebesgue measure, respectively, on the following product spaces:
$$\Sim\Omega\mathstrut^s_r,\quad \Sim\Omega_r. \eqno (3)$$
Finally, let 
$${\cal L}^s_r,\quad \buildrel\wedge\over{\cal L}\mathstrut^s_r,\quad \buildrel\vee\over{\cal L}\mathstrut^s_r,\quad 
\buildrel\wedge\over{\cal L}_r,\quad \buildrel\sim\over{\cal L}\mathstrut^s_r,\quad \buildrel\sim\over{\cal L}_r$$
be the respective $\sigma$-algebras of Lebesgue measurable subsets of the product spaces given in (2) and (3).
Let $(X,{\cal A},\nu)$ denote any one of the following measure spaces.
$$(\Omega^s_r,{\cal L}^s_r,m^s_r),\;(\Wedge\Omega\mathstrut^s_r,\buildrel\wedge\over{\cal L}\mathstrut^s_r,\Wedge m\mathstrut^s_r),
\;(\Vee\Omega\mathstrut^s_r,\buildrel\vee\over{\cal L}\mathstrut^s_r,\Vee m\mathstrut^s_r),\;
(\Wedge\Omega_r,\buildrel\wedge\over{\cal L}_r,\Wedge m_r),\;(\Sim\Omega\mathstrut^s_r,\buildrel\sim\over{\cal L}\mathstrut^s_r,\Sim m\mathstrut^s_r),
\;(\Sim\Omega_r,\buildrel\sim\over{\cal L}_r,\Sim m_r).$$
Then $(X,{\cal A},\nu)$ satisfies the hypothesis of Theorem 3.7, and hence in Definition 3.1 we may take $(X,{\cal M}_\mu,\mu)=(X,{\cal A},\nu)$. The 
measures in (1) are probability measures, and hence the we may apply Definition 3.1 to the following functions.
$$\varphi^{s,b}_{r,a}:\Omega^s_r\to\Lambda^{s,b}_{r,a},\;\varphi^s_{r,a}:\,\Wedge\Omega\mathstrut^s_r\to\Lambda^s_{r,a},\;
\varphi^{s,b}_r:\,\Vee\Omega\mathstrut^s_r\to\Lambda^{s,b}_r,\;\varphi_{r,a}:\;\Wedge\Omega_r\to\Lambda_{r,a}.$$
Therefore we may use these functions to construct, respectively, the image measure spaces:
$$(\Lambda^{s,b}_{r,a},{\cal M}^s_r,\mu^{s,b}_{r,a}),\;(\Lambda^s_{r,a},\buildrel{\;\wedge}\over{\cal M}\mathstrut^s_r,\mu^s_{r,a}),\;
(\Lambda^{s,b}_r,\buildrel{\;\vee}\over{\cal M}\mathstrut^s_r,\mu^{s,b}_r),\;(\Lambda_{r,a},\buildrel{\;\wedge}\over{\cal M}\mathstrut_r,\mu_{r,a}).$$
According to Theorem 3.5, the functions
$$\varphi^s_r:\,\Sim\Omega\mathstrut^s_r\to\Lambda^s_r,\quad \varphi_r:\,\Sim\Omega_r\to\Lambda_r$$
satisfy condition (a) of Definition 3.1, and hence we may use these functions to construct, respectively, the measure spaces:
$$(\Lambda^s_r,\buildrel{\;\sim}\over{\cal M}\mathstrut^s_r,\mu^s_r),\;(\Lambda_r,\buildrel{\;\sim}\over{\cal M}\mathstrut_r,\mu_r).$$
This completes the definition of the measure spaces
$$\displaylines{(\Lambda^{s,b}_{r,a},{\cal M}^s_r,\mu^{s,b}_{r,a}),\;(\Lambda^s_{r,a},\buildrel{\;\wedge}\over{\cal M}\mathstrut^s_r,\mu^s_{r,a}),\;
(\Lambda^{s,b}_r,\buildrel{\;\vee}\over{\cal M}\mathstrut^s_r,\mu^{s,b}_r),\cr
(\Lambda_{r,a},\buildrel{\;\wedge}\over{\cal M}\mathstrut_r,\mu_{r,a}),\;
(\Lambda^s_r,\buildrel{\;\sim}\over{\cal M}\mathstrut^s_r,\mu^s_r),\;(\Lambda_r,\buildrel{\;\sim}\over{\cal M}\mathstrut_r,\mu_r).\cr}$$
\medskip

\medskip
\medskip
\medskip
\centerline{\bf 4. UNIFORM MEASURE ON $\Lambda^{\,\square}_{\,\square}$}
\medskip
\medskip
\medskip
\noindent
In this section we select the functions
$$\lambda^{s,b}_{r,a}:[0,1]\to I^{s,b}_{r,a},\quad \lambda^s_{r,a}:[0,1]\to I^s_{r,a},\quad 
\lambda_r:\hbox{\b R}\to \hbox{\b R}$$
to be, respectively, the unique affine mappings of $[0,1]$ onto $I^{s,b}_{r,a}$, $I^s_{r,a}$, and the identity mapping
on {\b R}. These mappings give rise to ``the uniform probability measure on'' $\Lambda^{\,\square}_{\,\square}$ and 
``Lebesgue measure'' on  $\Lambda^{\,\square}_{\,\square}$.

\proclaim Definition 4.1. Let $0\le r<s$ be arbitrary. Let $(a,b)$ be any pair in $L^s_r$. For $\xi\in [0,1]$, define 
$$\eqalign{\lambda^{s,b}_{r,a}(\xi)&=\cases{[(a-b)+c(s-r)]\xi+b-{1\over 2}c(s-r),& $a\le b$;\cr
                           \mathstrut\cr
                                  [(b-a)+c(s-r)]\xi+a-{1\over 2}c(s-r),& $b\le a$,\cr}\cr
				   \mathstrut\cr
           \lambda^s_{r,a}(\xi)&=2c(s-r)\xi+a-c(s-r).\cr}$$
Then the affine function 
$$\lambda^{s,b}_{r,a}:[0,1]\to I^{s,b}_{r,a}$$
satisfies the conditions of Definition 2.3; and the affine function
$$\lambda^s_{r,a}:[0,1]\to I^s_{r,a}$$
satisfies the conditions of Definition 2.8. Hence we may use the functions $\lambda^{s,b}_{r,a}$ and $\lambda^s_{r,a}$ to construct the mappings
$$\varphi^{s,b}_{r,a}:\Omega^s_r\to\Lambda^{s,b}_{r,a},\;\varphi^s_{r,a}:\,\Wedge\Omega\mathstrut^s_r\to\Lambda^s_{r,a},\;
\varphi^{s,b}_r:\,\Vee\Omega\mathstrut^s_r\to\Lambda^{s,b}_r,\;\varphi_{r,a}:\;\Wedge\Omega_r\to\Lambda_{r,a}.$$
By Definition 3.2, these mappings give rise to the following image measure spaces.
$$(\Lambda^{s,b}_{r,a},{\cal M}^s_r,\mu^{s,b}_{r,a}),\;(\Lambda^s_{r,a},\buildrel{\;\wedge}\over{\cal M}\mathstrut^s_r,\mu^s_{r,a}),\;
(\Lambda^{s,b}_r,\buildrel{\;\vee}\over{\cal M}\mathstrut^s_r,\mu^{s,b}_r),
\;(\Lambda_{r,a},\buildrel{\;\wedge}\over{\cal M}\mathstrut_r,\mu_{r,a}). \eqno (1)$$
We call each measure $\lambda^{\,\square}_{\,\square}$ in (1) {\bf the uniform probability measure on} $\Lambda^{\,\square}_{\,\square}$.
Define 
$$\lambda_r:\hbox{\b R}\to \hbox{\b R}$$
by $\lambda_r(\xi)=\xi$, $\xi\in\hbox{\b R}$. Then $\lambda_r$ satisfies the conditions of Definition 2.10.
Hence, according to Definition 3.2, the functions
$$\lambda^{s,b}_{r,a}:[0,1]\to I^{s,b}_{r,a},\; \lambda^s_{r,a}:[0,1]\to I^s_{r,a},\; \lambda_r:\hbox{\b R}\to\hbox{\b R}$$
give rise to the following image measure spaces:
$$(\Lambda^s_r,\buildrel{\;\sim}\over{\cal M}\mathstrut^s_r,\mu^s_r),\;(\Lambda_r,\buildrel{\;\sim}\over{\cal M}\mathstrut_r,\mu_r).\eqno (2)$$
The measures  $\lambda^{\,\square}_{\,\square}$ in (2) are called {\bf Lebesgue measure on} 
$\Lambda^{\,\square}_{\,\square}$.

\vskip .3cm
\centerline{\bf References}
\vskip .2cm
\parindent=1cm
\item{\bf [HS]} Hewitt E., and Stromberg K. {\it Real and Abstract Analysis}, second edition, Springer Verlag, New York,
1969.
\item{\bf [M]} Munkres, J.R. {\it Topology, a first course}, Prentice Hall, Englewood Cliffs, New Jersey, 1975.
\bigskip

Department of Mathematics, University of Iowa, Iowa City, Iowa 52242

{\it E-mail address}: {\bf baker@math.uiowa.edu}

\vfill

\end